\documentclass{article}
\usepackage{amssymb, amsmath, amsthm, amsfonts,xcolor, enumerate, fullpage}
\usepackage{tikz, float}
\usetikzlibrary{positioning}

\usepackage[backend=biber,style=ieee,sorting=nyt]{biblatex}
\usepackage{cases}

\def\JJ{\mathbb{J}}

\def\phi{\varphi}

\def\AK{\color{red}}
\def\LZ{\color{blue}}
\def\ol{\overline}
\bibliography{references}
\begin{document}
\newpage 
\title{\bf {Tur\'{a}n number of four vertex-disjoint cliques}}

\author{Alexandr Kostochka\thanks{
Department of Mathematics, University of Illinois at Urbana-Champaign, Urbana, IL 61801, USA.
{\tt kostochk@illinois.edu.} The research of this author is partially supported by NSF RTG grant DMS-1937241.
}
\and {Dadong Peng}\thanks{
Department of Mathematics, University of Illinois at Urbana-Champaign, Urbana, IL 61801, USA.
{\tt dadongp2@illinois.edu.}}
\and Liang Zhang\thanks{Center for Combinatorics and LPMC, Nankai University, Tianjin 300071, P.R.China.  E-mail: {\tt 1120240010@mail.nankai.edu.cn}. }}
\date{}
\maketitle  

\begin{center}
\parbox{0.9\hsize}
{\small {\bf Abstract.}\ \ Given a graph $H$, the {\it Tur\'{a}n number} ${\rm ex}(n,H)$ of $H$ is the maximum number of edges of an
$n$-vertex simple graph containing no $H$ as a subgraph. Let
$kK_p$ denote the disjoint union of $k$ copies of the complete graph $K_p$. In this paper, utilizing the idea of the proof of the Hajnal–Szemer\'{e}di Theorem and discharging, we determine the values of ${\rm ex}(n,4K_p)$ for all $n\geq 4p+1$ and $p\ge 3$ and describe all extremal examples. Note that the result for $n<4p$ is trivial and for $n=4p$ was known before.
}
\end{center}

\begin{flushright} Dedicated to the 85th Birthday of Gyula Katona\end{flushright}

\section{Introduction}

\hskip\parindent  Graphs in this paper are finite, undirected and simple. Terms and notation not defined here are from \cite{BondyMurty1976}. We use $|S|$ to denote the cardinality of $S$. The vertex set and edge set of a graph $G$ are denoted by $V(G)$ and $E(G)$, respectively. The number of edges of a graph $G$ is denoted by $e(G)$.   For  a graph $G$, $v\in V(G)$ and $H\subseteq G$(respectively, $S\subseteq V(G)$), the set of neighbors of $v$ in $H$(respectively, $S$) is denoted by $N_H(v)$(respectively, $N_S(v)$).  We call  $d_G(v)=|N_G(v)|$   {\em the degree of $v$ in $G$}.  For vertex subsets $V_1$ and $V_2$ of a graph $G$, we let $E[V_1,V_2]$ denote the set of edges of $G$ with one endvertex in $V_1$ and the other  in $V_2$, and  let $e[V_1,V_2]=|E[V_1,V_2]|$. We denote by $\delta(G)$ and $\Delta(G)$ the minimum degree and the maximum degree of a graph $G$. 

Let $\ol{G}$ denote the complement of a graph $G$. The {\it independence number} of a graph $G$ is denoted by $\alpha(G)$. For a graph $G$ and $V'\subseteq V(G)$, the subgraph of $G$ induced by $V'$ is denoted by $G[V']$. We use $G\cup H$ to denote the disjoint union of graphs $G$ and $H$. We denote by $G\vee H$ the join operation of graphs $G$ and $H$, which means adding  edges between each vertex in $G$ and each vertex in $H$. We use $K_p$ to denote the complete graph on $p$ vertices, $\ol{K}_p$ to denote the edgeless graph on $p$ vertices and  $S_{\ell}$ to denote a star on $\ell+1$ vertices. For given graphs $G$ and $H$, $G$ is {\it $H$-free} if $G$ does not contain  $H$ as a subgraph. 

The Tur\'{a}n number of a graph $H$, denoted by ${\rm ex}(n,H)$, is the maximum number of edges in an $H$-free graph on $n$ vertices. An $n$-vertex graph with ${\rm ex}(n,H)$ edges not containing a copy of $H$ is an {\em extremal} graph for $H$. 
Let ${\rm EX}(n,H)$ denote the set of $n$-vertex extremal graphs for $H$. 

Let $T_{n, p}$ denote the complete balanced $p$-partite graph $K_{n_1, \ldots, n_p}$, where $n_1+\ldots +n_p=n$ and $\lfloor\frac{n}{p}\rfloor\le n_i\le\lceil\frac{n}{p}\rceil$ for $1\le i\le p$. Let $t_{n, p}$ be the number of edges of $T_{n, p}$. The famous result of Tur\'{a}n \cite{Turan1941} is Theorem 1.1. 

{\bf Theorem 1.1.} \cite{Turan1941}\ \ {\it $${\rm ex}(n,K_p)=t_{n,p-1}
\;\mbox{and \; ${\rm EX}(n,K_p)=\{T_{n, p-1}\}$.}
$$ }

We will   use the following notation: $\ol{\rm ex}(n,G):={n\choose 2}- {\rm ex}(n,G)$, $\ol{t}_{n,p}:={n\choose 2}-t_{n,p}$ and
$\ol{\rm EX}(n,H):=\{\ol{G}: G\in {\rm { EX} }(n,H)\}$. {This is because most of our arguments will be for the complements of the original graphs. It will be convenient to find four disjoint independent $p$-sets in the complement instead of looking for four cliques of size $p$ in the original graph. }

In these terms, Theorem 1.1  has the following corollary:

{\bf Corollary 1.1.1.}\ \ {\it Let $G$ be a graph on $n$ vertices with $\alpha(G)\le p-1$. Then $e(G)\ge \ol{t}_{n,p-1}$ and $G=\ol{T}_{n,p-1}$ if $e(G)=\ol{t}_{n,p-1}$.}

This extends the Mantel Theorem \cite{Mantel1907}, which shows ${\rm ex}(n,K_3)=\lfloor\frac{n^2}{4}\rfloor$. Tur\'{a}n numbers and extremal graphs for other graphs are widely studied in extremal graph theory. Gyula Katona made important contributions to the studies of 
Tur\'{a}n numbers of graphs and the applications of them to probability and geometry, see e.g. \cite{Katona1985,Katona1983,Katona2000,KatonaXiaoXiaoZamora2022,KatonaXiao2024}.

There are only a few graphs whose Tur\'{a}n number is determined exactly, see \cite{ErdosGallai1959,FaudreeSchelp1975,FurediGunderson2015,Kopylov1977}.


In 1959, Erd\H{o}s and Gallai \cite{ErdosGallai1959} determined the value of ${\rm ex}(n,kK_2)$ for all positive $n$ and $k$.

{\bf Theorem 1.2.} \cite{ErdosGallai1959}\ \ {\it Let $k\ge 1$ and $n\ge 2k$.
   Then $${\rm ex}(n,kK_2)=\max\left\{{2k-1\choose 2}, {k-1\choose 2}+(k-1)(n-k+1)\right\}$$ and ${\rm EX}(n,kK_2)\subseteq \{K_{2k-1}\cup\ol{K}_{n-2k+1},K_{k-1}\vee\ol{K}_{n-k+1}\}$.


} 

Determining the Tur\'{a}n number of vertex-disjoint copies of cliques was studied by Erdős \cite{Erdos1962}. Some years later, Moon \cite{Moon1968} and Simonovits \cite{Simonovits1968} determined ${\rm ex}(n,kK_p)$ for sufficiently large $n$.

{\bf Theorem 1.3.} \cite{Moon1968} \cite{ Simonovits1968}\ \ {\it For each fixed $k, p$ and sufficiently large $n$, $${\rm ex}(n,kK_p)={k-1\choose 2}+t_{n-k+1,p-1}+(k-1)(n-k+1), \;\mbox{and\; ${\rm EX}(n,kK_p)=\{K_{k-1}\vee T_{n-k+1,p-1}\}$. }$$
}


The remaining question is to determine the value of ${\rm ex}(n, kK_p)$ for every $n$ and $k$. There are only a few cases when the Tur\'{a}n number ${\rm ex}(n, kK_p)$ is known exactly for $n\ge kp$. In 2022, Chen, Lu and Yuan \cite{ChenLuYuan2022} determined the Turán number of two vertex-disjoint copies of $K_p$ completely. 

{\bf Theorem 1.4.} \cite{ChenLuYuan2022}\ \ {\it If $p\ge 3$, then 
$\ol{\rm ex}(n,2K_p)=\left\{
\begin{array}{ll} 3(n-2p+1), & \mbox{if $2p\le n\le 3p-2$,}\\ {\ol{t}_{n-1,p-1}},  & \mbox{if $n\ge 3p-1$.}\end{array}\right.$}

Zhang and Yin \cite{ZhangYin2024Klambda} determined the value of ${\rm ex}(n,K_p\cup K_q)$ for all $n,q$ and $p=2,3$. Later, Hu \cite{Hu2024} determined ${\rm ex}(n,K_p\cup K_q)$ for all $p$ and $q$. Furthermore, Luo \cite{Luo2025} determined all extremal graphs for $K_p\cup K_q$.

{\bf Theorem 1.5.} \cite{Hu2024}\ \ {\it Let $n$, $p$, $q$ be positive integers with $q>p\ge 3$ and $n\ge p+q$. Then $$\ol{\rm ex}(n,K_p\cup K_q)=\left\{\begin{array}{lll} 3(n-p-q+1), & \mbox{if $n\le p+q+\max\{2p-q,\lfloor\frac{p}{2}\rfloor-1$\},}\\
\ol{t}_{n,q-1},  & \mbox{if $n>p+q+\max\{2p-q,\lfloor\frac{p}{2}\rfloor-1$\}.}\end{array}\right.$$}

Brualdi and Mellendorf \cite{BrualdiMellendorf1994} and independently Zhang \cite{Zhang2024DisjointCliques} determined ${\rm ex}(kp,kK_p)$ for all $k$ and $p$.

{\bf Theorem 1.6.} \cite{BrualdiMellendorf1994} \cite{Zhang2024DisjointCliques}\ \ {\it Let $p\ge 3$, $k\ge 1$ and $H$ be an extremal graph for $kK_p$ on $kp$ vertices. Then
$$\ol{\rm ex}(kp,kK_p)=\left\{\begin{array}{ll}
{{k+1}\choose2}, & \mbox{if $k\le 2p-2$,}\\
kp-p+1, & \mbox{if $k\ge 2p-1$.}\end{array}\right.$$ Moreover, either $\overline{H}=K_{k+1}\cup\overline{K}_{kp-k-1}$ for $k\le2p-2$, or $\overline{H}\in\{K_{1,x}\cup(kp-p-x+1)K_2\cup(2p-kp+x-3)K_1:kp-2p+3\le x\le kp-p+1\}$ for $k\ge2p-2$.
}

Recently, Zhang and Yin \cite{ZhangYin2024ThreeCliques} and independently Zhang \cite{Zhang2024DisjointCliques} determined the value of ${\rm ex}(n,3K_p)$ for all $n$. Yan, Zhang and Peng~\cite{Yan2024DisjointCliques} also obtained this result and described the extremal graphs.

{\bf Theorem 1.7.} \cite{Yan2024DisjointCliques}\cite{Zhang2024DisjointCliques}\cite{ZhangYin2024ThreeCliques}\ \ {\it If $p\ge 3$, then
$$\ol{\rm ex}(n,3K_p)=\left\{\begin{array}{ll}
6, & \mbox{if $n=3p$,}\\
5(n-3p+1), & \mbox{if $3p+1\le n\le 5p-2$,}\\
\ol{t}_{n-2,p-1}, & \mbox{if $n\ge
5p-1$.}\end{array}\right.$$}

In this paper, we further determine ${\rm ex}(n,4K_p)$ and ${\rm EX}(n,4K_p)$ for all $n\ge 4p$. Let 
\begin{equation}\label{Jnxy}
 \mbox{ { $J_n(x,y)=xK_7\cup yK_8\cup\ol{K}_{n-7x-8y}$ }, where $x,y$ are non-negative integers with $7x+8y\le n$,   }
\end{equation}
and
\begin{equation}\label{Jnxyp}
 { \JJ_n(p)=\left\{{J_n(x,y)}: 3x+4y=n-4p+1\ \mbox{ and }\ x+y\le p-1\right\}.}
\end{equation}
Moreover, for $7p-6\le n\le 8p-8$, we let \begin{equation}\label{Jn'}
{\JJ'_n(p)=\JJ_n(p)\cup \left\{(8p-8-n)K_7\cup (n-7p+6)K_8\cup S_7\right\}.}
\end{equation}

{ Our main result is:}

{\bf Theorem 1.8.} \ \ {\it If $p\ge 3$, then 

1.\ $\ol{\rm ex}\left(4p+1,4 K_{p}\right)=15$, and 
$\ol{\rm EX}(4p+1,4K_p)=\{K_6\cup\ol{K}_{4p-5}\}$;

2.\ $\ol{\rm ex}\left(4p+2,4 K_{p}\right)=21$,
and 
$\ol{\rm EX}(4p+2,4K_p)=\{K_7\cup\ol{K}_{4p-5}\}$;

3.\ $\ol{\rm ex}\left(4p+3,4 K_{p}\right)=28$, 
$\ol{\rm EX}(4p+3,4K_p)=\{K_8\cup\ol{K}_{4p-5}\}$ when $p\geq 4$ and
{ $\ol{\rm EX}(15,4K_3)=\{K_8\cup\ol{K}_{7}, K_7\cup S_7\}$};

4.\ $\ol{\rm ex}\left(4p+4,4 K_{p}\right)=  36$, $\ol{\rm EX}(4p+4,4K_p)=\{K_6\cup K_7\cup\ol{K}_{4p-9}, K_9\cup\ol{K}_{4p-5}\}$ when $p\geq 4$,  { $\ol{\rm ex}\left(16,4 K_{3}\right)=35$ and
$\ol{\rm EX}(16,4K_3)=\{ K_8\cup S_7\}$ };

5.\ $\ol{\rm ex}(n,4K_p)=\left\{\begin{array}{lll} 
7(n-4p+1), & \mbox{if $4p+5\le n\le 7p-2$,}\\
\ol{t}_{n-3,p-1},  & \mbox{if $n\ge 7p-1$,}\end{array}\right.$

and
{ $$\ol{\rm EX}(n,4K_p)=\left\{\begin{array}{lll} 
\JJ_n(p), & \mbox{if $4p+5\le n\le 7p-7$,}\\
\JJ'_n(p), & \mbox{if $7p-6\le n\le 8p-8$,}\\
\left\{\ol{K}_3\cup \ol{T}_{n-3,p-1}\right\},  & \mbox{if $n\ge 8p-7$.}\end{array}\right.$$}



}

For some results on the Tur\'{a}n number for disjoint copies of other graphs, we refer the readers to \cite{BielakKieliszek2014_3P4,BielakKieliszek2016_2P5,BushawKettle2011,CamposLopes2018,Gorgol2011,LanLiShiTu2019,LiYinLi2022,LidickyLiuPalmer2013,YuanZhang2017,YuanZhang2021}.

The structure of the paper is as follows. In the next section, we show the sharpness of our bound, cite the known results and lemmas that we will use, and state Theorem 2.7, a slight variation of an important case of our main theorem. In Section 3, we prove Theorem 2.7, and in Section 4, we use it to prove Theorem 1.8. We conclude the paper with some remarks in Section 5.

\section{Preliminaries}

\hskip\parindent 
The idea of the proof of Theorem 1.8 is to consider the problem in the complementary form: we look at graphs that contain no four disjoint independent $p$-sets. Let $s=n-4p+1$. First we discuss constructions providing upper bounds on $\ol{\rm ex}(n,4K_p)$.

{\bf Claim 2.1.}\ \ {\it Every graph in  $\JJ_n(p)$ or $\JJ'_n(p)$ contains no four disjoint independent $p$-sets and has $7s$ edges.}

{\bf Proof of Claim 2.1.}\ \ {Let $J\in \JJ_n(p)$, say $J=J_n(x,y)$ with $x+y\leq p-1$ and $3x+4y=s$.
Since $n=4p-1+s$, every independent set of size $p$ in $J$ contains at least $p-x-y$ vertices from $(n-7x-8y)K_1$. As $4(p-x-y)-(n-7x-8y)=4p+3x+4y-n$ and $3x+4y=n-4p+1$, $4(p-x-y)-(n-7x-8y)=1$. Thus $4(p-x-y)>(n-7x-8y)$. So $J$ has no four disjoint independent $p$-sets. Also $e(J)=21x+28y=7(3x+4y)=7(n-4p+1)=7s$. 

{ For $7p-6\le n\le 8p-8$, we consider $(8p-8-n)K_7\cup (n-7p+6)K_8\cup S_7$. Since $S_7$ has no four disjoint independent $2$-sets, $(8p-8-n)K_7\cup (n-7p+6)K_8\cup S_7$ contains no four disjoint independent $p$-sets. Moreover, $e((8p-8-n)K_7\cup (n-7p+6)K_8\cup S_7)=7(n-4p)+7=7(n-4p+1)=7s$. }\ $\Box$ }

By Claim 2.1,  
$\ol{\rm ex}(n,4K_p)\le 7s$ when $4p+5\leq n\leq 8p-8$. 
Similarly, every independent set of size $p$ in $\ol{K}_3\cup \ol{T}_{n-3,p-1}$  contains at least one vertex from $\ol{K}_3$, and hence  $\ol{K}_3\cup \ol{T}_{n-3,p-1}$ has no four disjoint independent $p$-sets.
In the same vein for $2\leq s\leq 5$, graph
$K_{s+4}\cup\ol{K}_{4p-5}$ does not contain four disjoint independent sets of size $p$. We leave it to the reader to check that
graphs $K_8\cup\ol{K}_{7}, K_7\cup S_7 $ and $ K_8\cup S_7$ have no four disjoint independent $3$-sets
and that graph $K_6\cup K_7\cup\ol{K}_{4p-9}$ has no four disjoint independent $p$-sets. 

An important part of our lower bounds is to prove that for $ 4p+5\leq n\leq 7p-2$ 
every $n$-vertex graph with at most $7s$ edges other than the  graphs described above has four disjoint independent sets of size $p$. 

To prove this, we need some lemmas and known results of equitable coloring. An {\it equitable $k$-coloring} of a graph $G$ is a proper $k$-coloring in which any two color classes differ in size by at most one. In 1970, Hajnal and Szemer\'{e}di \cite{HajnalSzemeredi1970} proved the following well-known result.

{\bf Theorem 2.2.} \cite{HajnalSzemeredi1970}\ \ {\it Every graph with maximum degree at most $r$ has an equitable $(r+1)$-coloring.}

Chen, Lih and Wu \cite{ChenLihWu1994}  proposed the following conjecture and confirmed it for $\Delta(G)\le 3$.

{\bf Conjecture 2.3.} \cite{ChenLihWu1994}\ \ {\it If $G$ is an $r$-colorable graph with $\Delta(G)\le r$, then either $G$ has an equitable $r$-coloring, or $r$ is odd and $K_{r,r}\subseteq G$.}

Kierstead and Kostochka \cite{KiersteadKostochka2012} confirmed the conjecture for $\Delta=4$.

{\bf Theorem 2.4.} \cite{KiersteadKostochka2012}\ \ {\it Let $r\le 4$ and $G$ be an $r$-colorable graph with $\Delta(G)\le r$. Then either $G$ has an equitable $r$-coloring or $r$ is odd and $G$ contains $K_{r,r}$.}

Based on Theorem 2.2 and the idea in \cite{ChenLuYuan2022}, Zhang \cite{Zhang2024DisjointCliques} obtained the following lemma.

{\bf Lemma 2.5.} \cite{Zhang2024DisjointCliques}\ \ {\it Let $G_p(n)=K_{k-1}\vee T_{n-k+1,p-1}$ and $n_0\ge kp$, where $p\ge 3$ and $k\ge 2$. If ${\rm ex}(n_0, kK_p)=e(G_p(n_0))$, then ${\rm ex}(n, kK_p)=e(G_p(n))$ for every $n\ge n_0$.}

Using Theorem 2.4, Zhang and Yin \cite{ZhangYin2024ThreeCliques} obtained the following result.

{\bf Lemma 2.6.} \cite{ZhangYin2024ThreeCliques}\ \ {\it Let $p\ge 3$, $1\le s\le 2p-1$ and $n=3p-1+s$. Then $$\ol{\rm ex}(n,3K_p)=\left\{\begin{array}{lll} 6, & \mbox{if $s=1$,}\\
5s,  & \mbox{if $2\le s \le 2p-1$.}\end{array}\right.$$}

In the spirit of this lemma, we prove the following theorem which helps to prove an important part of  our main result, Theorem 1.8.

{ {\bf Theorem 2.7.} \ \ {\it For $p\ge 3$, let $G$ be a graph on $n=4p-1+s$ vertices with $1\le s\le 3p-1$  and $\Delta(G)\le 6$. 
(a) If $s\equiv 0\pmod 3$ and $|E(G)|\le 7s$, then $G$ contains four disjoint independent $p$-sets or $G=\frac{s}{3}K_7\cup\ol{K}_{n-\frac{7s}{3}}$;\\
(b) If $s\not\equiv 0\pmod 3$ and $|E(G)|\le 7s+1$, then $G$ contains four disjoint independent $p$-sets or $s\equiv 2\pmod 3$ and $G=\frac{s-2}{3}K_7\cup K_6\cup\ol{K}_{n-\frac{7s+4}{3}}$;\\
(c) If $s=3$ and $e(G)\le 22$, then $G$ contains four disjoint independent $p$-sets or $K_7\subseteq G$.} }



After proving Theorem 2.7, we consider several special cases for small $n$, and use induction and case analysis together with Theorem 2.7 to complete the proof of Theorem 1.8.

\section{ Proof of Theorem 2.7}

In the first subsection, we introduce an auxiliary digraph whose vertices are vertex subsets of a graph $G$, and derive its basic properties when $G$ is a counter-example for Theorem 2.7. There will be four possible cases for a structure of such a digraph. In each of the four subsequent subsections, we consider one such case, and 
show that $G$ actually is not a counter-example to Theorem 2.7.

\subsection{An auxiliary digraph and its basic properties}
\hskip\parindent We need some definitions. Given a partition $V_1,\ldots,V_k$ of $V(G)$, define an auxiliary digraph $D$ with vertices $V_1,\ldots,V_k$, so that $V_iV_j(1\le i,j\le k)$ is a directed edge if and only if some vertex $x\in V_i$ has no neighbors in $V_j$. In this case, we say that $x$ is {\it movable} to $V_j$. 

Let $V_{i_1},V_{i_2},\ldots, V_{i_s}\in V(D)$ and $v_{i_\ell}\in V_{i_\ell} $ for each $\ell\in \{1,2,\ldots,s-1\}$. 
The sequence $V_{i_1}V_{i_2}\ldots V_{i_s}$ is called an {\em accessible path from $V_{i_1}$ to $V_{i_s}$} if we can find a vertex shifting as above, where for each $\ell\in \{1,2,\ldots,s-1\}$, $v_{i_\ell}$ is movable to $V_{i_{\ell+1}}$. 
We also say that $V_j$ is inaccessible for $V_i$ if there is no accessible path  from $V_j$ to $V_i$.

Set a target set $V_k$. Call $V_i\in V(D)$ {\it accessible} if there is an accessible path from $V_i$ to $V_k$ in $D$. Note that $V_k$ is trivially accessible. Let $\mathcal{A}$ be the set of accessible classes, $\mathcal{B}= V(D)-\mathcal{A}$, $A=\bigcup_{V_i\in \mathcal{A}}V_i$ and
$B=\bigcup_{V_j\in \mathcal{B}}V_j$. 

For $v_i\in V_i$ and $v_j\in V_j$, the edge $v_iv_j$ between two classes is called a {\em solo edge} for $v_i$ if it is the only edge from $v_i$ to the class $V_j$. In this case, $v_j$ is called a 
{\em solo neighbor} of $v_i$. 


{\bf Claim 3.1.}\ \ {\it For each $V_i\in \mathcal{B}$ and each $V_j\in \mathcal{A}$, every $v_i\in V_i$ has a neighbor in $V_j$.}

{\bf Proof of Claim 3.1.}\ If there are $V_i\in \mathcal{B}$, $V_j\in \mathcal{A}$ and $v_i\in V_i$ such that $N_{V_j}(v_i)=\emptyset$, then $V_iV_j\in E(D)$. Since $V_j\in \mathcal{A}$, there is a directed path $V_j\ldots V_k$ in $D$. Thus $V_iV_j\ldots V_k$ is a directed path in $D$. This means that $V_i\in \mathcal{A}$, a contradiction. \ $\Box$




Let $G$ be a counterexample to Theorem 2.7 with the minimum number of edges {over all $p\geq 3$ and $1\leq s\leq 3p-1$}: it satisfies all conditions without containing four disjoint independent $p$-sets. In particular, this implies that for any $uv\in E(G)$, $G-uv$ contains four disjoint independent $p$-sets, and one of these sets contains both $u$ and $v$. When we move $u$ out of the union of these sets, $V(G)$ can be divided into 5 classes $V_1,V_2,V_3,V_4,V_5$, such that $V_2,V_3,V_4,V_5$ are four independent sets with $|V_2|=|V_3|=|V_4|=p$, $|V_5|=p-1$  and $V_1$ is the set of remaining vertices. Notice that $|V_1|=s$.

In the auxiliary digraph $D$, we let class $V_5$ be our target set, let $\mathcal{A}$ be composed of accessible sets and $\mathcal{B}$ be composed of inaccessible sets. By definition, $V_5\in \mathcal{A}$. 
Class $V_1$ is always in $\mathcal{B}$, since otherwise there is a directed path from $V_1$ to $V_5$. This directed path provides a vertex shifting which makes $V_2,V_3,V_4,V_5$ be $4$ independent sets of size $p$.

{\bf Claim 3.2.}\ \ {\it If $v\in V_i\in\mathcal{A}\setminus\{V_5\}$ is a solo neighbor of a vertex $x\in V_1$ and $D$ contains a directed path from $V_j$ to $V_5$ that avoids $V_i$, then $N_{V_j}(v)\not=\emptyset$.
{ In particular, always $N_{V_1}(v)\not=\emptyset$.} }


{\bf Proof of Claim 3.2.}\ \  If $N_{V_j}(v)=\emptyset$, then $v$ is movable to $V_j$. We can move $x$ to $V_i$ and move $v$ to $V_j$. The directed path from $V_j$ to $V_5$  avoiding $V_i$ in $D$ provides a vertex shifting from $V_j$ to $V_5$ that does not involve  $V_i$.  By Claim 3.1, it also avoids $V_1$. 
Hence $G$ contains four disjoint independent $p$-sets, a contradiction.\ $\Box$

{\bf Claim 3.3.}\ \ {\it No vertex $v\in V_i\in\mathcal{A}$ is a solo neighbor of two non-adjacent vertices in $V_1$.}

{\bf Proof of Claim 3.3.}\ \ Assume that $v\in V_i\in\mathcal{A}$ is a solo neighbor of two non-adjacent vertices $v_1,v_1'\in V_1$. Since $V_i\in\mathcal{A}$, there is a directed path $P=V_i\ldots V_5$, where $V(P)\subseteq\mathcal{A}$.
This directed path $P$ provides a vertex shifting from $V_i$ to $V_5$. 
If $i\neq 5$, let $v'\in V_i$ be movable to the successor of $V_i$ in $P$. 
By Claim 3.2, $v\neq v'$. We do the shifting along $P$, then move $v_1,v_1'$ to $V_i$ and move $v$ to $V_1$. Now $G$ contains four disjoint independent $p$-sets, a contradiction.\ $\Box$

{\bf Claim 3.4.}\ \ {\it If $\Delta(G[V_1])\le 2$, then $1\le s \le 3(p-1)$ and $e(G[V_1])\ge\ol{t}_{s,p-1}$.}

{\bf Proof of Claim 3.4.}\ \ Suppose $\Delta(G[V_1])\le 2$. Then $G[V_1]$ is $3$-colorable. Thus if $|V_1|=s\geq 3(p-1)+1$, then $G[V_1]$ contains an independent $p$-set. In this case, $G[V_1\cup V_2\cup V_3\cup V_4]$ contains four disjoint independent $p$-sets, a contradiction. Thus $1\le s \le 3(p-1)$. Note that $\alpha(G[V_1])\le p-1$, otherwise $G[V_1]$ contains an independent $p$-set, and so $G$ contains four disjoint independent $p$-sets, a contradiction. By Corollary 1.1.1, $e(G[V_1])\ge\ol{t}_{s,p-1}$. \ $\Box$


\begin{equation}\label{|B|}
   \mbox{Among all  partitions as above, choose one with the smallest  $|\mathcal{B}|$.} 
\end{equation}

In the next four subsections, we consider the four different possibilities for $|\mathcal{B}|$.

\subsection{\it $|\mathcal{B}|=1$.}

\hskip\parindent By Claim 3.1, $e[V_1,V_2\cup V_3\cup V_4\cup V_5]\ge 4s$. Since $V_1\in \mathcal{B}$, $\Delta(G)\le 6$ and every vertex in $V_1$ has a neighbor in each class in $\mathcal{A}$, $\Delta(G[V_1])\le 2$. We first prove some properties of the solo neighbors of vertices in $V_1$, and then will use discharging to prove Theorem 2.7 in this case.

{\bf Claim 3.5.}\ \ {\it We can move vertices between classes in $\mathcal{A}$ so that $|\mathcal{B}|=1$ and at least two classes in $\mathcal{A}\setminus\{V_5\}$ are inneighbors of $V_5$.}

{\bf Proof of Claim 3.5.}\ \
Assume that there is exactly one class in $\mathcal{A}\setminus\{V_5\}$ that is an inneighbor of $V_5$, say $V_2V_5\in E(D)$, $V_3V_5\not\in E(D)$ and $V_4V_5\not\in E(D)$. This implies that any directed path from $V_3$ or $V_4$ to $V_5$ must go through $V_2$. {In particular, for $i\in \{3,4\}$, $D$ contains a path $P_i$
from $V_i$ to $V_2$.}
Without loss of generality, we assume $V_3V_2\in E(D)$. Then we move a movable vertex {$v_2\in V_2$} to $V_5$, and let {$V_2-v_2$} be a new destination set. 
{In the new digraph $D'$, paths $P_3$ and $P_4$  still exist, and $v_2$ witnesses an edge from $V_5+v_2$ to $V_2-v_2$.
So, the number of accessible sets does not decrease, and $V_3$ and $V_5+v_2$ are inneighbors of $V_2-v_2$ in $D'$.}
\ $\Box$

{\bf Claim 3.6.}\ \ {\it Every vertex in $V_i\in\mathcal{A}\setminus\{V_5\}$ that is the solo neighbor of a vertex in $V_1$ has neighbors in each class in $\mathcal{A}\setminus\{V_i\}$.}

{\bf Proof of Claim 3.6.}\ \ Let $i\in\{2,3,4\}$, and $v_i\in V_i$ be a solo neighbor of $v_1\in V_1$. By Claim 3.2, $v_i$ has neighbors in $V_5$. By Claim 3.5, 
it is enough to consider the following two cases.

{\bf Case 1.}\ \ {\it Each class in $\mathcal{A}\setminus\{V_5\}$ is an inneighbor of $V_5$.}
Since $V_2V_5, V_3V_5, V_4V_5\in E(D)$, we are done by Claim 3.2.

{\bf Case 2.}\ \ {\it There are exactly two classes in $\mathcal{A}\setminus\{V_5\}$ that are inneighbors of $V_5$.}
We may let $V_2V_5\not\in E(D)$ and $V_3V_5, V_4V_5\in E(D)$. Since $V_2\in\mathcal{A}$, without loss of generality, we assume that $V_2V_3\in E(D)$. If $V_4V_3\in E(D)$, let $v_3'\in V_3$ with $N_{V_5}(v_3')=\emptyset$. Then we set $V'_3=V_5\cup\{v_3'\}$ and $V'_5=V_3\setminus\{v_3'\}$.
Now $V_2V'_5, V_4V'_5, V'_3V'_5\in E(D)$, and we have Case 1. 
 So, assume that $V_4V_3\not\in E(D)$.

For any $v_2\in V_2$, by $V_3V_5, V_4V_5\in E(D)$ and Claim 3.2, $N_{V_3}(v_2)\not=\emptyset$ and $N_{V_4}(v_2)\not=\emptyset$. For any $v_4\in V_4$, by $V_3V_5\in E(D)$, $V_2V_3V_5\subseteq D$ and Claim 3.2, $N_{V_2}(v_4)\not=\emptyset$ and $N_{V_3}(v_4)\not=\emptyset$. Next, we consider $v_3\in V_3$. Since $V_4V_5\in E(D)$, by Claim 3.2, $N_{V_4}(v_3)\not=\emptyset$. 

Then we only need to prove that $N_{V_2}(v_3)\not=\emptyset$ for each $v_3\in V_3$. If $V_2V_4\in E(D)$, then $V_2V_4V_5\subseteq D$, and by Claim 3.2, we have $N_{V_2}(v_3)\not=\emptyset$. So we assume that $V_2V_4\not\in E(D)$. Then $e[V_2,V_4]\ge p$. Recall that $V_2V_5, V_4V_3\not\in E(D)$, so $e[V_2,V_5]\ge p$ and $e[V_4,V_3]\ge p$. 

By Claim 3.1, every vertex in $V_1$ has at least one neighbor in each class in $\mathcal{A}$. Assume that there are $x$ vertices in $V_1$ that have a solo neighbor in $V_2$, and $s-x$ vertices in $V_1$ that have at least two neighbors in $V_2$. Recall that $\Delta(G[V_1])\le 2$. It follows from Theorem 2.2 and Claim 3.3 that there are at least $\lceil\frac{x}{3}\rceil$ vertices in $V_1$ that have distinct solo neighbors in $V_2$. Since $N_{V_3}(v_2)\not=\emptyset$, 
$$e[V_1,V_2]+ e[V_2,V_3]\ge x+2(s-x)+\lceil\frac{x}{3}\rceil=2s-\lfloor\frac{2x}{3}\rfloor\ge 2s-\frac{2x}{3}.$$ 
So, since
$x\le s$, $e[V_1,V_2]+ e[V_2,V_3]\ge\frac{4s}{3}$. Similarly, it follows from $N_{V_5}(v_3)\not=\emptyset$ and $N_{V_5}(v_4)\not=\emptyset$ that  $e[V_1,V_3]+ e[V_3,V_5]\ge\frac{4s}{3}$ and $e[V_1,V_4]+ e[V_4,V_5]\ge\frac{4s}{3}$. Moreover, since $\Delta(G[V_1])\le 2$, by Claim 3.4, $1\le s\le 3(p-1)$ and $e(G[V_1])\ge\ol{t}_{s,p-1}$. Thus, we have the following:

$\begin{cases}
{ e(G[V_1])\ge\ol{t}_{s,p-1}, }\\
e[V_1,V_2]+ e[V_2,V_3]\ge\frac{4s}{3},\quad
e[V_1,V_3]+ e[V_3,V_5]\ge\frac{4s}{3},\quad
e[V_1,V_4]+ e[V_4,V_5]\ge\frac{4s}{3},\\
e[V_1\cup V_2,V_5]\ge s+p, \quad
 e[V_4,V_3]\ge p,\quad e[V_2,V_4]\ge p.
\end{cases}$

\noindent Summing  these inequalities we get $e(G)\ge 5s+3p+\ol{t}_{s,p-1}$.

Since 
{ $\ol{t}_{s,p-1}=s-(p-1)$ for $p\le s<2(p-1)$ and $\ol{t}_{s,p-1}=2s-3(p-1)$ 
for $2(p-1)\le s\le 3(p-1)$, } { we have}
$$e(G)\ge5s+3p+{s\choose 2}-t_{s,p-1} =\left\{\begin{array}{lll}{ 5s+3p>7s+1 }, & \mbox{if $0\le s<p$,}\\ 5s+3p+s-(p-1)={6s+2p+1>7s+1 }, & \mbox{if $p\le s<2(p-1),$}\\ 5s+3p+2s-3(p-1)={ 7s+3>7s+1 }, & \mbox{if $2(p-1)\le s\le 3(p-1)$.}\end{array}\right.$$ This contradiction to $e(G)\leq 7s+1$ yields $N_{V_2}(v_3)\not=\emptyset$.\ $\Box$

\medskip
Consider the following discharging procedure. At the start, each $v\in V(G)$ has charge $ch(v)=0$ and each $e\in E(G)$ has charge $ch(e)=1$. So, $\sum_{x\in V(G)\cup E(G)}ch(x)=|E(G)|\leq  7s+1$. Now we will move the charges between edges and vertices without changing the total sum as follows.

Step 1.\ \ Every edge $xy\in E(G)$ such that  $x\in V_1,y\not \in V_1$ gives charge $1$ to $x$.

Step 2.\ \ Every edge $xy\in E(G)$ such that  $x,y\in V_1$ gives charge 1/2 to $x$ and 1/2 to $y$.

Step 3.\ \ { Every edge $xy\in E(G)$ such that  $x\in A\setminus V_5$ 
 and $y\in V_5$ 
gives charge $1$ to $x$.}

Step 4.\ \ Every edge $xy\in E(G)$ such that $x,y\in A\setminus V_5$ 
gives charge 1/2 to $x$ and 1/2 to $y$.

Step 5.\ \ Every vertex $v$ in $V_i\in\mathcal{A}$ distributes  its charge equally between the vertices in $V_1$ for which $v$ is the solo neighbor  in $V_i$.

Let the charge of each $x\in V(G)\cup E(G)$ after Step $j$ be denoted by $ch_j(x)$. Then $ch_4(e)= 0$ for every $e\in E(G)$ and $ch_5(v)\geq 0$ for every $v\in A$. 
Therefore 
\begin{equation}\label{charge}
    \parbox{14.5cm}{\begin{center}$e(G)\ge\sum_{v\in V_1}ch_5(v),$\end{center}
    
    and  equality holds if and only if for each edge $xy$ with $\{x,y\}\subseteq A$, each of $x$ and $y$ either is in $V_5$ or is a solo neighbor of a vertex in $V_1$.
    }
\end{equation}

For $0\leq j\leq 2$, let $V_{1,j}=\{v\in V_1: d_{V_1}(v)=j\}$.

{\bf Claim 3.7.}\ \ {\it For each $v\in V_1$, 

1)\ If $v\in V_{1,0}$, then $ch_5(v)\ge 8$. Equality holds  only if 
$v$ has exactly one solo neighbor in $A\setminus V_5$. 

2)\ If  $v\in V_{1,1}$, then $ch_5(v)\ge 7.5$. Equality holds  only if all solo neighbors of $v$ in $A\setminus V_5$ are also solo neighbors of its neighbor in $V_1$.

3)\ If  $v\in V_{1,2}$, then $ch_5(v)\ge 7$. Equality holds only if the solo neighbors of $v$  in $V_2, V_3, V_4$ are also solo neighbors of its two neighbors in $V_1$.

4)\ If  $v\in V_{1,2}$, then its two neighbors in $V_1$ are adjacent to each other. Moreover, $ch_5(v)=7$ for all $v\in V_{1,2}$.}

{\bf Proof of Claim 3.7.}\ \ Recall that $\Delta(G[V_1])\le 2$. Thus, we consider the following cases.

{\bf Case 1.}\ \ {\it $v\in V_{1,0}$.} By Claim 3.1, $v$ has  neighbors
in each of $V_2, V_3, V_4, V_5$. By Claim 3.3, solo neighbors of $v$ are not solo neighbors of other vertices in $V_1$. Assume that $v$ has exactly $x$  solo neighbors in $A\setminus V_5$. Since $\Delta(G)\le 6$,  $x\ge 1$. Then $ch_1(v)\geq x+(3-x)\cdot 2+1$. Let $V_i\in\mathcal{A}\setminus\{V_5\}$ and $v_i\in V_i$ be a solo neighbor of $v$. By Claim 3.6, $v_i$ has neighbors in each class in $\mathcal{A}\setminus\{V_i\}$. Hence $ch_4(v_i)\geq \frac{1}{2}+\frac{1}{2}+1=2$. On Step 5, each solo neighbor will send all its charge to $v$. Therefore, $$ch_5(v)\ge x\cdot(1+2)+(3-x)\cdot 2+1=7+x\ge 8,$$ { and equality holds only if $x=1$. 
 }

{\bf Case 2.}\ \ {\it $v\in V_{1,1}$.}
Let $N_{V_1}(v)=\{v'\}$. By Claim 3.3,  each solo neighbor of $v$ in $A$ can be a solo neighbor only of $v$ and $v'$. Since $d_G(v)-d_{V_1}(v)\le 6-1=5$, $v$ has solo neighbors in at least $3$ classes of $\mathcal{A}$.

If $v$ has $4$ solo neighbors in ${A}$, then $|N(v)|=5$. In this case, let $N(v)=\{v', a_2, a_3, a_4, a_5\}$, where $a_i\in V_i$ for $2\le i\le 5$. Then $ch_1(v)=4$. At Step 2, $v$ receives charge $\frac{1}{2}$ from  edge $vv'$. By Claim 3.6,  $ch_4(a_i)\ge\frac{1}{2}+\frac{1}{2}+1=2$ for $2\le i\le 4$. At Step 5,  $a_i$ ($2\le i\le 4$)  sends at least half of its charge to $v$. Therefore, $$ch_5(v)\ge 4+\frac{1}{2}+3\cdot\frac{1}{2}\cdot 2=7.5>7.$$

If $v$ has $3$ solo neighbors, then $v$ has $5$ neighbors in $A$. Without loss of generality, we may assume that  $a_2'\in V_2$  and $a_3'\in V_3$ are solo neighbors of $v$. Then $ch_1(v)=5$. At Step 2, $v$ receives charge $\frac{1}{2}$ from  edge $vv'$. By Claim 3.6, $ch_4(a_i')\ge\frac{1}{2}+\frac{1}{2}+1=2$ for $2\le i\le 3$. 
At Step 5, $a_i'$ ($2\le i\le 3$) sends at least half of its charge to $v$. Therefore, 
$$ch_5(v)\ge 5+\frac{1}{2}+2\cdot\frac{1}{2}\cdot 2=7.5.$$ Note that $ch_5(v)=7.5$ 
 only if $v$ and $v'$ have the same solo neighbors in $A\setminus V_5$.
 

{\bf Case 3.}\ \ {\it $v\in V_{1,2}$.}
Let $N_{V_1}(v)=\{a, b\}$. For $2\le i\le 5$, let $v_i$ be a neighbor of $v$ in $V_i\in\mathcal{A}$. Since $\Delta(G)\le 6$ and $d_{V_1}(v)=2$, $v_i$ is the solo neighbor of $v$. Hence
$ch_1(v)=4$. At Step 2, $v$ receives charge $\frac{1}{2}$ from each of the  edges $va$ and $vb$.
By Claim 3.6, $ch_4(v_i)\ge\frac{1}{2}+\frac{1}{2}+1=2$ for $2\le i\le 4$. By Claim 3.3, $v_i$ ($2\le i\le 4$) can be a solo neighbor of only $v, a$ and $b$. At Step 5, each $v_i$ ($2\le i\le 4$) sends at least $\frac{1}{3}$ of its charge to $v$. Therefore, $$ch_5(v)\ge 4+1+3\cdot\frac{1}{3}\cdot 2=7,$$ 
and $ch_5(v)=7$ 
only if each $v_i$ ($2\le i\le 4$) is the common solo neighbor of $v, a$ and $b$. This proves 3).

Furthermore, by 3), $ch_1(v)=4$ and $ch_4(v_i)\ge\frac{1}{2}+\frac{1}{2}+1=2$ for $2\le i\le 4$. If $ab\not\in E(G)$, then by Claim 3.3, $v_i$ is not a common solo neighbor of $a$ and $b$ for all $2\le i\le 4$. Then $ch_5(v)\ge 4+1+3\cdot\frac{1}{2}\cdot 2=8.$ By 2)--3), $ch_5(v)+ch_5(a)+ch_5(b)\geq 8+7.5+7.5=23.$ Also by 1)-3), since $e(G)\geq \sum_{v\in V_1}ch_5(v)$, { $$e(G)\geq 7(|V_{1}\setminus \{v,a,b\}|)+23 >7s+1,$$ a contradiction. Thus,
$ab\in E(G)$.}

If $ch_5(v)\not=7$, then by 3), at least one vertex in $\{v_2,v_3,v_4\}$ is not a common solo neighbor of $a$ and $b$. This implies that some vertex in $\{v, a, b\}$ has a solo neighbor in $V_2\cup V_3\cup V_4$ that is not a solo neighbor of the other two. Then $ch_5(v)+ch_5(a)+ch_5(b)\ge 2\cdot(4+1+2\cdot\frac{1}{3}\cdot 2+\frac{1}{2}\cdot 2)+(4+1+2\cdot\frac{1}{3}\cdot2+2)=23.$ By 1)-3), since $e(G)\geq \sum_{v\in V_1}ch_5(v)$, $$e(G)\ge 8\cdot|V_{1,0}|+7.5\cdot|V_{1,1}|+23+7\cdot(|V_{1,2}|-3)\ge 7s+2,$$ a contradiction.$\qquad\Box$

{\bf Claim 3.8.}\ \ {\it { If $e(G)\leq 7s+1$,}
then\\
(a) $|V_{1,0}|+\frac{1}{2}|V_{1,1}|\leq 1$, and 
$|V_{1,2}|=3\lfloor s/3\rfloor$;\\
(b) 
$ch_5(v)=
\left\{\begin{array}{ll}
 8    & \mbox{for each } v\in V_{1,0},\\
      7.5   & \mbox{for each } v\in V_{1,1},\\
      7   & \mbox{for each } v\in V_{1,2}; 
\end{array}\right.$\\  
(c) $G[V_{1,2}]=\frac{|V_{1,2}|}{3}K_3$, and for every $v\in V_{1,2}$, $G[\{v\}\cup N_{G\setminus V_5}(v)]=K_6$  and each vertex of this $K_6$ has a neighbor in $V_5$;\\
(d) every $v\in V_{1,2}$ is in a $K_7$.
}

{\bf Proof of Claim 3.8.}\ \ By part 4) of Claim 3.7,  every $v\in V_{1,2}$  has two neighbors in $V_1$ that are adjacent to each other.
So $|V_{1,2}|$ is divisible by 3. Also by Claim 3.7, 
\begin{equation}\label{012}
  e(G)\geq \sum_{v\in V_1}ch_5(v)\geq 8|V_{1,0}|+7.5|V_{1,1}|+7|V_{1,2}|=7s+|V_{1,0}|+\frac{1}{2}|V_{1,1}|.
\end{equation}
As $e(G)\leq 7s+1$ and $|V_{1,2}|$ is divisible by 3, this implies (a). 

{ By Claim 3.7, $ch_5(v)=7$ for each $v\in V_{1,2}$. If $s\equiv 0\pmod 3$, then by (a), $V_{1,0}=V_{1,1}=\emptyset$. If $s\not\equiv 0\pmod 3$, then by (a), $|V_{1,0}|+\frac{1}{2}|V_{1,1}|=1$. By~\eqref{012}, this implies (b).}

Let $v\in V_{1,2}$, $N_{V_1}(v)=\{v',v''\}$ and let
$v_i$ be a neighbor of $v$ in $V_i$ for $2\leq i\leq 5$. Denote
$W(v):=\{v, v', v'', v_2, v_3, v_4\}$.
Assume that $G[\{v_2, v_3, v_4\}]\not=K_3$, { say $v_2v_3\notin E(G)$.} By Claims 3.6 and 3.7(3), since $\Delta(G)\le 6$, we may assume $N_{V_3}(v_2)=\{u_3\}$ and $u_3\not=v_3$. We switch $v'$ with $v_2$, and get $\{v',v_3,v_4,v_5, u_3\}\subseteq N_{A}(\{v, v_2, v''\})$. By Claim 3.7, $ch_5(v)>7$,  contradicting  (b). Thus $G[\{v_2, v_3, v_4\}]=K_3$. This implies $G[W(v)]=K_6$.
Furthermore, each $x\in W(v)$ has a neighbor in $V_5$: for $x\in \{v,v',v''\}$ this holds
since $V_1=B$, and for $x\in \{v_2,v_3,v_4\}$ this holds by Claim 3.6.
This proves (c).

Suppose now that $W(v)$ has at least two neighbors in $V_5$. By (a), $s\geq 3$.
By the symmetry between vertices in $W(v)$, we may assume that $v_5,v'_5\in V_5$ and $vv_5,v'v_5'\in E(G)$. Let $G'$ be obtained from $G\setminus W(v)$ by adding edge $v_5v'_5$. Then 
 \begin{equation}\label{10}
\mbox{(i) $|E(G')|= |E(G)|-20$ and (ii) $\Delta(G')\leq 6$.}
\end{equation}
 
Suppose $G'$ contains four disjoint independent $(p-1)$-sets, say $W_1,\ldots,W_4$. If there are four distinct vertices
$w_1,\ldots,w_4\in W(v)$ such that for $j=1,2,3,4$,  $w_j$ has no neighbors in $W_j$, then the sets
$W'_j=W_j+w_j$ form four disjoint independent $p$-sets, a contradiction. Otherwise,  
{consider the auxiliary bipartite graph $F$ with parts $W(v)$ and $U=\{1,2,3,4\}$, where we have an edge from a $w\in W(v)$ to a $j\in U$ if and only if $w$ has no neighbor in $W_j$.
Since  each $w\in W(v)$ has at most one neighbor outside of $W(v)$ in $G$, by Hall's Theorem for $F$,}
there is a $j\in [4]$ such that $N(W_j)\supseteq W(v)$.
But then $v_5,v'_5\in W_j$, contradicting the fact that $W_j$ is independent.

So, $G'$ does not contain four disjoint independent $(p-1)$-sets.
It has {$n'=n-6=4(p-1)-1+(s-2)$} vertices. Since $s\geq 3$, $s':=s-2\geq 1$
and by~\eqref{10}(ii), $\Delta(G')\leq 6$.
So, if $p\geq 4,$
then by the minimality of $G$, $e(G')\geq 7(s-2)$. 
And if $p=3$, then applying Theorem~1.2 for $k=4$ we again get
$e(G')\geq 7(s-2)$.
Thus, in both cases
by~\eqref{10}(i), $e(G)\geq 7(s-2)+20>7s+1$, a contradiction.
This proves (d).
$\Box$

\smallskip
We now finish the proof of the case $|\mathcal{B}|=1$. If 
$s\equiv 0\pmod 3$, then by Claim 3.8(a), $V_{1,0}=V_{1,1}=\emptyset$ and 
by Claim 3.8(d), all vertices in $V_1$ are in $K_7$, three vertices in each $K_7$. {If $e(G)\le 7s$, then by~\eqref{012}, $e(G)=\sum_{v\in V_1}ch_5(v)=7s$.} Now by~\eqref{charge}
, $G$ has no other edges, thus the theorem holds in this case. {If $e(G)=7s+1$, then by Claim 3.8(d), $K_7\subseteq G$.}

{ If $s\not\equiv 0\pmod 3$, then by Claim 3.8(a), $|V_{1,0}|+\frac{1}{2}|V_{1,1}|=1$. Since $e(G)\le 7s+1$, by~\eqref{012}, $e(G)=\sum_{v\in V_1}ch_5(v)=7s+1$.} Suppose now $s\equiv 1\pmod 3$. By Claim 3.8(a), $V_{1,1}=\emptyset$ and $|V_{1,0}|=1$, say $V_{1,0}=\{y\}$. As we observed before, since $\Delta(G)\leq 6$, $y$ has a solo neighbor in some $V_i$ for $2\leq i\leq 4$, say, $v_2$ is a solo neighbor of $y$ in $V_2$. By Claim 3.6, $v_2$ has a neighbor in $V_3$, say $v_3$.
Then since $e(G)=\sum_{v\in V_1}ch_5(v)=7s+1$, by Claim 3.8(b), $v_3$ is also a solo neighbor of a vertex in $V_1$, say $v_1$. But by Claim 3.8(d), $v_1\notin V_{1,2}$. Thus
$v_1=y$, and this contradicts part 1) of Claim 3.7.

Finally, suppose $s\equiv 2\pmod 3$. Now by Claim 3.8(a), $V_{1,0}=\emptyset$ and $|V_{1,1}|=2$, say $V_{1,1}=\{y,y'\}$.
Again, $y$ has a solo neighbor in some $V_i$ for $2\leq i\leq 4$, say, $v_2$ is a solo neighbor of $y$ in $V_2$. Again by Claim 3.6, $v_2$ has  neighbors in $V_3$ and $V_4$, say $v_3\in V_3$ and $v_4\in V_4$.

Since $e(G)=\sum_{v\in V_1}ch_5(v)=7s+1$, by Claim 3.8(b), 
$v_3$ and $v_4$ are also solo neighbors of some vertices in $V_1$, say $v_3$ is a solo neighbor of $v_1$ and $v_4$ is a solo neighbor of $v_1'$. By Claim 3.8(d), $v_1\notin V_{1,2}$.
So $v_1\in \{y,y'\}$, and by part 2) of Claim 3.7, the solo neighbors of $y$ and $y'$ are the same. Using a similar argument for $v_4$ and symmetry, we obtain that $G[\{y,y',v_2,v_3,v_4\}]=K_5$. Let $W=\{y,y',v_2,v_3,v_4\}$. Each $v\in W$ has a neighbor in $V_5$: $y$ and $y'$ because they are in $B$, and $v_2,v_3,v_4$ because they are solo neighbors of $y$. Thus, we already have $21\frac{s-2}{3}$ edges in $K_7$s and $15$ edges incident to $W$; together, $7(s-2)+15=7s+1$ edges.
So, $G$ has no other edges, and each $w\in W$ has exactly one neighbor in $V_5$. Since
$W$ has exactly one vertex in each of $V_2,V_3,V_4$, every copy of $K_7$ has exactly one vertex in each of $V_2,\ldots,V_5$ and $|V_5|<p$, for each $i\in \{2,3,4\}$ the set $V_i$ contains an isolated vertex, say $z_i$.
 Let $N_{V_5}(y)=\{v_5\}$ and $N_{V_5}(y')=\{v'_5\}$. If $v_5$ has no neighbors in $V_i$ for some $i\in \{2,3,4\}$, then after moving $y$ and $z_i$ to $V_5$ and  $v_5$ to $V_i-z_i$  we get 
four disjoint independent $p$-sets, a contradiction. Thus $v_5$ has a neighbor $x_i\in V_i$ for each  $i\in \{2,3,4\}$. Since this edge must be one of the $15$ incident to $W$,
$x_i=v_i$. Repeating this argument with $v'_5$ in place of $v_5$ we get that $v_5'=v_5$ and
$G[W\cup \{v_5\}]=K_6$. It follows that $G=K_6\cup\frac{s-2}{3}K_7\cup \ol{K}_{n-\frac{7s+4}{3}}$, as claimed.

\subsection{\it $|\mathcal{B}|=2$}

\hskip\parindent We may let $\mathcal{B}=\{V_1,V_2\}$. By Claim 3.1, $e[V_1,V_3\cup V_4\cup V_5]\ge 3s$ and $e[V_2,V_3\cup V_4\cup V_5]\ge 3p$. If $V_3V_5\notin E(D)$, then since $V_3\in \mathcal{A}$, $V_3V_4,V_4V_5\in E(D)$. So, 
there exists $v_3\in V_3$ that can
be moved to $V_4$ and $v_4\in V_4$ that can be moved to $V_5$.
 We move $v_4$ to $V_5$ to get the new target set $V_4-v_4$ and sets $V_3$ and $V_5+v_4$ that have vertices movable to $V_4-v_4$. Thus, we can assume that $V_3V_5,V_4V_5\in E(D)$. 

{\bf CASE 1.}\ \ {\it $V_1V_2\notin E(D)$.} Since $V_1V_3, V_1V_4, V_1V_5\notin E(D)$, $\Delta(G[V_1])\le 2$. Consider the following discharging procedure. At the start, each $v\in V(G)\setminus V_2$ has charge $ch(v)=0$ and each $e\in E[V_1, V_3\cup V_4\cup V_5]\cup E(G[V_3\cup V_4\cup V_5])$ has charge $ch(e)=1$. So, $\sum_{x\in V(G)\setminus V_2\cup E[V_1, V_3\cup V_4\cup V_5]\cup E(G[V_3\cup V_4\cup V_5])}ch(x)=e[V_1, V_3\cup V_4\cup V_5]+e(G[V_3\cup V_4\cup V_5])$. Now we will move the charges between edges and vertices without changing the total sum as follows.

Step 1.\ \ Every edge $xy\in E[V_1, V_3\cup V_4\cup V_5]$ such that $x\in V_1, y\not\in V_1$ gives charge $1$ to $x$.

{ Step 2.\ \ Every edge $xy\in E(G[V_3\cup V_4\cup V_5])$ such that  $x\in A\setminus V_5$ and $y\in V_5$ gives charge $1$ to $x$.}


{ Step 3.\ \ Every edge $xy\in E(G[V_3\cup V_4\cup V_5])$ such that $x,y\in A\setminus V_5$ gives charge 1/2 to $x$ and 1/2 to $y$.}


Step 4.\ \ Every vertex $v$ in $V_i\in\mathcal{A}$ distributes its charge equally between the vertices in $V_1$ for which $v$ is the solo neighbor  in $V_i$.

Let the charge of each $x\in V(G)\setminus V_2\cup E[V_1, V_3\cup V_4\cup V_5]\cup E(G[V_3\cup V_4\cup V_5])$ after Step $j$ be denoted by $ch_j(x)$. Then { $ch_3(e)=0$} for every $e\in E[V_1, V_3\cup V_4\cup V_5]\cup E(G[V_3\cup V_4\cup V_5])$ and $ch_4(v)\geq 0$ for every $v\in A$. Therefore $$e[V_1, V_3\cup V_4\cup V_5]+e(G[V_3\cup V_4\cup V_5])\ge\sum_{v\in V_1}ch_4(v),$$ and the equality holds if and only if for each edge $xy$ of $E(G[V_3\cup V_4\cup V_5])$, each of $x$ and $y$ either is in $V_5$ or is a solo neighbor of a vertex in $V_1$.

{\bf Claim 3.9.}\ \ {\it For each $v\in V_1$, $ch_4(v)\ge 4$.}

{\bf Proof of Claim 3.9.}\ \ Recall that $\Delta(G[V_1])\le 2$. Then we have the following cases.

{\bf Case {\romannumeral 1.}}\ \ {\it $d_{V_1}(v)=0$.} By Claim 3.1, $v$ has neighbors
in each of $V_3, V_4, V_5$. By Claim 3.3, those solo neighbors cannot be the solo neighbors of another vertex in $V_1$. Assume that $v$ has exactly $x$ solo neighbors in $A\setminus V_5$. Since $\Delta(G)\le 6$, 
$ch_1(v)\geq x+(2-x)\cdot 2+1$. 

Let $V_i\in\{V_3, V_4\}$ and $v_i\in V_i$ be a solo neighbor of $v$. By $V_3V_5, V_4V_5\in E(D)$ and Claim 3.2, $v_i$ has neighbors in each class in $\{V_3,V_4,V_5\}\setminus\{V_i\}$. Hence $ch_3(v_i)\geq\frac{1}{2}+1$. On Step 4, each solo neighbor will send all its charge to $v$. Therefore, $$ch_4(v)\ge x\cdot(\frac{1}{2}+2)+(2-x)\cdot 2+1=5+\frac{x}{2}>4.$$

{\bf Case {\romannumeral 2.}}\ \ {\it $d_{V_1}(v)=1$.} Let $N_{V_1}(v)=\{v'\}$. By Claim 3.3,  each solo neighbor of $v$ in $A$ can be a solo neighbor only of $v$ and $v'$. Since $d_G(v)-d_{V_1}(v)\le 6-1=5$ and $V_1V_2\not\in E(D)$, $v$ has solo neighbors in at least two classes of $\mathcal{A}$.

If $v$ has $3$ solo neighbors in $A$, then $|N_A(v)|=3$. Thus $ch_1(v)=3$. In this case, let $N_{A}(v)=\{a_3, a_4, a_5\}$, where $a_i\in V_i$ for $3\le i\le 5$. By $V_3V_5, V_4V_5\in E(D)$ and Claim 3.2, $ch_3(a_i)\ge\frac{1}{2}+1=\frac{3}{2}$ for $3\le i\le 4$. At Step 4,  $a_i$ ($3\le i\le 4$) sends at least half of its charge to $v$. Therefore, $$ch_4(v)\ge 3+2\cdot\frac{1}{2}\cdot\frac{3}{2}=4.5>4.$$ If $v$ has two solo neighbors, then $|N_A(v)|\ge 4$. Thus $ch_4(v)\ge ch_1(v)\ge 4$.

{\bf Case {\romannumeral 3.}}\ \ {\it $d_{V_1}(v)=2$.} Let $N_{V_1}(v)=\{a, b\}$. For $3\le i\le 5$, let $v_i$ be a neighbor of $v$ in $V_i\in\mathcal{A}$. Since $\Delta(G)\le 6$, $d_{V_1}(v)=2$ and $V_1V_2\not\in E(D)$, $v_i$ is the solo neighbor of $v$. Hence
$ch_1(v)=3$. By $V_3V_5, V_4V_5\in E(D)$ and Claim 3.2, $ch_3(v_i)\ge\frac{1}{2}+1=\frac{3}{2}$ for $3\le i\le 4$. By Claim 3.3, $v_i$ ($3\le i\le 4$) can be a solo neighbor of only $v, a$ and $b$. At Step 4, each $v_i$ ($3\le i\le 4$) sends at least $\frac{1}{3}$ of its charge to $v$. Therefore, $$ch_4(v)\ge 3+2\cdot\frac{1}{3}\cdot\frac{3}{2}=4.\qquad \qquad \Box$$ 

By Claim 3.9, $e[V_1, V_3\cup V_4\cup V_5]+e(G[V_3\cup V_4\cup V_5])\ge 4s$. Since $\Delta(G[V_1])\le 2$, by Claim 3.4, $1\le s\le 3(p-1)$ and $e(G[V_1])\ge {s\choose 2}-t_{s,p-1}$. Since $V_1V_2\not\in E(D)$, $e[V_1,V_2]\ge s$. In summary, we have the following:

$\begin{cases}
e(G[V_1])\ge {s\choose 2}-t_{s,p-1},\qquad
e[V_1,V_2]\ge s,\\
e[V_1, V_3\cup V_4\cup V_5]+e(G[V_3\cup V_4\cup V_5])\ge 4s, \qquad
e[V_2, V_3\cup V_4\cup V_5]\ge 3p.
\end{cases}$

\noindent Then { \begin{equation}\label{B21}{e(G)\ge5s+3p+{s\choose 2}-t_{s,p-1}\geq\left\{\begin{array}{lll}5s+3p>7s+1, & \mbox{if $0\le s<p$,}\\ 6s+2p+1>7s+1, & \mbox{if $p\le s<2(p-1),$}\\ 7s+3>7s+1, & \mbox{if $2(p-1)\le s\le 3(p-1)$,}\end{array}\right.}\end{equation} }
a contradiction. 

{\bf CASE 2.}\ \ {\it $V_1V_2\in E(D)$.} 
Then there is $v_1\in V_1$  movable to $V_2$.
Since $\mathcal{B}=\{V_1,V_2\}$,  $\Delta(G[V_1\cup V_2])\le 3$. 

{\bf Claim 3.10.}\ \ {\it { If for some $i=3,4$ a vertex $v_i\in V_i$ is} a solo neighbor of $v\in V_1\cup V_2$, then $v_i$ has neighbors in each class of $\{V_3,V_4,V_5\}\setminus\{V_i\}$}. 

{\bf Proof of Claim 3.10.}\ \ If $v\in V_2$, then we move $v_1$ to $V_2$ and move $v$ to $V_1$. Thus we may assume $v\in V_1$. As $V_3V_5, V_4V_5\in E(D)$, by Claim 3.2, { we are done. } $\Box$

{\bf Claim 3.11.}\ \ {\it { No vertex in $ V_i\in\mathcal{A}\setminus\{V_5\}$  is } a solo neighbor of two non-adjacent vertices in $V_1\cup V_2$. 
  }

{\bf Proof of Claim 3.11.}\ \ Assume that $v\in V_i\in\mathcal{A}\setminus\{V_5\}$ is a solo neighbor of two non-adjacent vertices $v_x,v_y\in V_1\cup V_2$. If $v_x,v_y\in V_1$, then  we are done by Claim 3.3. Thus we may assume $v_y\in V_2$.


Recall that $V_iV_5\in E(D)$ for $i=3,4$. Let $v'\in V_i$ be movable to $V_5$. By Claim 3.10, $v'\not=v$. Then 
{ we move $v_x$ to $V_i$ and $v$ to $V_1$.}
After this, $v_y$ is movable to $V_i$ and $V_jV_5\in E(D)$ still holds for $j=3,4$. If $v_x\in V_1$, then $|\mathcal{B}|\le 1$, contradicting~\eqref{|B|}. If $v_x\in V_2$, then we move $v_1$ to $V_2$. Again $|\mathcal{B}|\le 1$. \ $\Box$


{ Denote
$E_0:=E[V_1\cup V_2, V_3\cup V_4\cup V_5]\cup E(G[V_3\cup V_4\cup V_5])$ and } consider the following discharging procedure. 
At the start, each $v\in V(G)$ has charge $ch(v)=0$ and each $e\in E_0$ has charge $ch(e)=1$. So, $$\sum_{x\in V(G)\cup E_0}ch(x)=e[V_1\cup V_2, V_3\cup V_4\cup V_5]+e(G[V_3\cup V_4\cup V_5]).$$ Now we will move the charges between edges and vertices without changing the total sum as follows.

Step 1.\ \ Every edge $xy\in E[V_1\cup V_2, V_3\cup V_4\cup V_5]$ such that $x\in V_1\cup V_2, y\not\in V_1\cup V_2$ gives charge $1$ to $x$.

{ Step 2.\ \ Every edge $xy\in E(G[V_3\cup V_4\cup V_5])$ such that $x\in A\setminus V_5$ and $y\in V_5$ gives charge $1$ to $x$.}


{ Step 3.\ \ Every edge $xy\in E(G[V_3\cup V_4\cup V_5])$ such that $x,y\in A\setminus V_5$ gives charge 1/2 to $x$ and 1/2 to $y$.}


Step 4.\ \ Every vertex $v$ in $V_i\in\mathcal{A}$ distributes its charge equally between the vertices in $V_1\cup V_2$ for which $v$ is the solo neighbor in $V_i$.

Let the charge of each $x\in V(G)\cup E_0$ after Step $j$ be denoted by $ch_j(x)$. Then $ch_3(e)=0$ for every $e\in E_0$ and $ch_4(v)\geq 0$ for every $v\in A$. Therefore 
$$ { |E_0|}\ge\sum_{v\in V_1\cup V_2}ch_4(v),$$ and the equality holds if and only if for each edge $xy$ of $E(G[V_3\cup V_4\cup V_5])$, each of $x$ and $y$ either is in $V_5$ or is a solo neighbor of a vertex in $V_1\cup V_2$.


{\bf Claim 3.12.}\ \ {\it For each $v\in V_1\cup V_2$, $ch_4(v)\ge\frac{15}{4}$.}

{\bf Proof of Claim 3.12.}\ \ By Claim 3.1, $v$ has neighbors in each of $V_3, V_4, V_5$. If $v$ has a solo neighbor $v_3\in V_3$ and a solo neighbor $v_4\in V_4$, then $|N_A(v)|\ge 3$. Thus $ch_1(v)\geq 3$. By Claim 3.10, $v_i$ has neighbors in each class in $\{V_3,V_4,V_5\}\setminus\{V_i\}$ for $i=3,4$. Hence $ch_3(v_i)\geq\frac{1}{2}+1=\frac{3}{2}$. Recall that $\Delta(G[V_1\cup V_2])\le 3$. By Claim 3.11, $v_i$ ($3\le i\le 4$) can be a solo neighbor of only $v$ and the vertices in $N_{V_1\cup V_2}(v)$. On Step 4, each solo neighbor will send at least $\frac{1}{4}$ of its charge to $v$. Therefore, $$ch_4(v)\ge 3+2\cdot\frac{1}{4}\cdot\frac{3}{2}=\frac{15}{4}.$$

Assume now that $v$ has the solo neighbor in exactly one class of $\{V_3, V_4\}$. We may let $v_3'\in V_3$ be a solo neighbor of $v$ and $v$ has at least two neighbors in $V_4$. Then $|N_A(v)|\ge 4$. Thus $ch_4(v)\ge ch_1(v)\geq 4>\frac{15}{4}$.

If $v$ has at least two neighbors in each of $V_3$ and $V_4$, then $|N_A(v)|\ge 5$. Thus $ch_4(v)\ge ch_1(v)\geq 5>\frac{15}{4}$. $\Box$

If $s<p$, then by Claim 3.12, 
{ \begin{equation}\label{B22}
e(G)\ge |E_0|\ge\frac{15}{4}(s+p)>7s+1,
\end{equation}}
a contradiction. Assume that $p\le s\le 3p-1$. If $G[V_1\cup V_2]$ contains two disjoint independent $p$-sets, then $G[V_1\cup V_2\cup V_3\cup V_4]$ contains four disjoint independent $p$-sets, a contradiction. If $G[V_1\cup V_2]$ contains no two disjoint independent $p$-sets, then by Theorem 1.4, $e(G[V_1\cup V_2])\ge 3(s-p+1)$. Thus
{ \begin{equation}\label{B23}
e(G)\ge e(G[V_1\cup V_2])+|E_0|
\ge 3(s-p+1)+\frac{15}{4}(s+p)={ \frac{27}{4}s+\frac{3}{4}p+3> 7s+1 },
\end{equation} }
a contradiction.



\subsection{\it $|\mathcal{B}|=3$.}

\hskip\parindent We may let $\mathcal{B}=\{V_1,V_2,V_3\}$. Then $\mathcal{A}=\{V_4,V_5\}$. This implies $V_4V_5\in E(D)$. By Claim 3.1, $e[V_1\cup V_2\cup V_3,V_4\cup V_5]\ge 2s+4p$.

{\bf CASE 1.}\ \ {\it $V_1V_2\in E(D)$ or $V_1V_3\in E(D)$.} We may assume $V_1V_2\in E(D)$. Let $v_1\in V_1$ be movable to $V_2$. 

{\bf Claim 3.13.}\ \ {\it No vertex in $ V_4$  is a solo neighbor of two vertices in $V_2$ 
.}

{\bf Proof of Claim 3.13.}\ \ Assume that $v\in V_4$ is a solo neighbor of  $v_2,v_2'\in V_2$. Then $v$ has a neighbor in $V_5$, since otherwise we can obtain four disjoint independent $p$-sets by moving $v_1$ to $V_2$, $v_2$ to $V_4$ and $v$ to $V_5$.
 Let $v'\in V_4$ be movable to $V_5$. Then $v'\not=v$. 
{We move $v_1$ to $V_2$, $v_2$ to $V_4$ and $v$ to $V_1$. Now}
$v_2'$ is movable to $V_4$ and $V_4V_5\in E(D)$ still holds. This makes $|\mathcal{B}|\le 2$, contradicting~\eqref{|B|}.\ $\Box$


Consider $e[V_2, V_4]+e[V_4, V_5]$. By Claim 3.1, every $v\in V_2$ has at least one neighbor in $V_4$. If $v$ has a solo neighbor $v_4\in V_4$, then $v_4$ has a neighbor in $V_5$. Otherwise we can obtain four disjoint independent $p$-sets by 
moving $v_1$ to $V_2$, $v$ to $V_4$ and $v_4$ to $V_5$.
 If there are $x$ vertices in $V_2$ having a solo neighbor in $V_4$, then $e[V_2, V_4]+e[V_4, V_5]\ge 2x+2(p-x)=2p$.

If $s<p$, then
{ \begin{equation}\label{B31}
e(G)\ge e[V_1, V_4\cup V_5]+(e[V_2, V_4]+e[V_4, V_5])+e[V_2, V_5]+e[V_3, V_4\cup V_5]\ge 2s+2p+p+2p=2s+5p> 7s+1,
\end{equation}}
a contradiction. Assume that $p\le s\le 3p-1$. If $G[V_1\cup V_2\cup V_3]$ contains three disjoint independent $p$-sets, then $G[V_1\cup V_2\cup V_3\cup V_4]$ contains four disjoint independent $p$-sets, a contradiction. If $G[V_1\cup V_2\cup V_3]$ contains no three disjoint independent $p$-sets, then by Theorem 1.7, $e(G[V_1\cup V_2\cup V_3])\ge 5(s-p+1)$. So, we have the following:

$\begin{cases}
e(G[V_1\cup V_2\cup V_3])\ge 5(s-p+1),\qquad
e[V_1, V_4\cup V_5]\ge 2s,\qquad e[V_3, V_4\cup V_5]\ge 2p,\\
e[V_2, V_4]+e[V_4, V_5]\ge 2p,\qquad e[V_2, V_5]\ge p.
\end{cases}$

\noindent Summing the inequalities we get
{ \begin{equation}\label{B32}
e(G)\ge 7s+5, 
\end{equation}}
a contradiction.

{\bf CASE 2.}\ \ {\it $V_1V_2,V_1V_3\notin E(D)$.} Then $\Delta(G[V_1])\le 2$. Let $v$ and $v'$ be two non-adjacent vertices in $V_1$. If $v_i\in V_2$ (or $V_3$) is a common solo neighbor of $v$ and $v'$, then we move $v$ to $V_2$ (or $V_3$) and move $v_i$ to $V_1$. Thus we obtain $V_1V_2\in E(D)$ (or $V_1V_3\in E(D)$) and accordingly turn to Case 1, which is done. Therefore, we can assume that any two non-adjacent vertices $v,v'\in V_1$ have no common solo neighbors in $V_2$ or $V_3$. 

Let $v_2\in V_2$ be a solo neighbor of $v\in V_1$. If $N_{V_3}(v_2)=\emptyset$, then we can switch $v$ and $v_2$ to turn into Case 1 (by $V_1V_3\in E(D)$), which is done. Thus we can assume that $v_2$ has a neighbor in $V_3$. Consider $e[V_1, V_2]+e[V_2, V_3]$. Since $V_1V_2\notin E(D)$, we may assume that there are $x$ vertices in $V_1$ that have a solo neighbor in $V_2$, and $s-x$ vertices in $V_1$ that have at least two neighbors in $V_2$. It follows from $\Delta(G[V_1])\le 2$ and Theorem 2.2 that there are at least $\lceil\frac{x}{3}\rceil$ vertices in $V_1$ that all have  distinct solo neighbors in $V_2$. Thus 
$$e[V_1,V_2]+ e[V_2,V_3]\ge x+2(s-x)+\lceil\frac{x}{3}\rceil=2s-\lfloor\frac{2x}{3}\rfloor\ge 2s-\frac{2x}{3}.$$ Since $x\le s$,  $e[V_1,V_2]+e[V_2,V_3]\ge\frac{4s}{3}$.

Now consider $e[V_1, V_4]+e[V_4, V_5]$. By Claim 3.3, no two non-adjacent vertices $v_1,v_1'\in V_1$ have a common solo neighbor in $V_4$. Let $v_4\in V_4$ be a solo neighbor of $v_1\in V_1$. By Claim 3.2, $v_4$ has neighbors in $V_5$. Similarly to the above discussion, $e[V_1, V_4]+e[V_4, V_5]\ge\frac{4s}{3}$. 

Since $\Delta(G[V_1])\le 2$, by Claim 3.4, $1\le s\le 3(p-1)$ and $e(G[V_1])\ge {s\choose 2}-t_{s,p-1}$. So, we have the following:

$\begin{cases}
e(G[V_1])\ge {s\choose 2}-t_{s,p-1},\qquad
e[V_1, V_2]+e[V_2, V_3]\ge\frac{4s}{3},\qquad e[V_1, V_3]\ge s,\\
e[V_1, V_4]+e[V_4, V_5]\ge\frac{4s}{3},\qquad e[V_1, V_5]\ge s,\qquad
e[V_2, V_4\cup V_5]\ge 2p, \qquad e[V_3, V_4\cup V_5]\ge 2p.
\end{cases}$

\noindent Then 
{ \begin{equation}\label{B33}
e(G)\ge\frac{14s}{3}+4p+{s\choose 2}-t_{s,p-1}=\left\{\begin{array}{lll}\frac{14s}{3}+4p>7s+1, & \mbox{if $0\le s< p$,}\\ \frac{17s}{3}+3p+1>7s+1, & \mbox{if $p\le s<2(p-1),$}\\ 
\frac{20s}{3}+p+3>7s+1, & \mbox{if $2(p-1)\le s\le 3(p-1)$,}\end{array}\right.
\end{equation}}
a contradiction. 


\subsection{\it $|\mathcal{B}|=4$}

\hskip\parindent In this case, $\mathcal{B}=\{V_1,V_2,V_3,V_4\}$. By Claim 3.1, $e[V_1,V_5]\ge s$ and $e[V_i,V_5]\ge p$ for $2\le i\le 4$. Define $\mathcal{B'}$ as
the set of all color classes $V_i$ in $\mathcal{B}$ 
such that there are no accessible paths from $V_1$ to  $V_i$.

\begin{equation}\label{min5}
 \mbox{Choose our partition with minimum $\sum\nolimits_{v\in V_5}d(v)$,
 and modulo this, with minimum 
 $|\mathcal{B'}|$.}
\end{equation}

{\bf Claim 3.14.}\ \ {\it 
If $v\in V_5$ is the solo neighbor of $u\in V_1\cup V_2\cup V_3\cup V_4$, then $d_G(u)\ge d_G(v)$.}

{\bf Proof of Claim 3.14.}\ \ Assume that $d_G(u)<d_G(v)$. It suffices to prove that after moving $u$ to $V_5$ and $v$ to $V_1$, we can switch vertices in ${B}$ so that $V_2, V_3$ and $V_4$ are independent sets of size $p$. 

If $u\in V_1$, then we move $u$ to $V_5$ and move $v$ to $V_1$. Since $d_G(u)<d_G(v)$, this contradicts~\eqref{min5}.
Assume that $u\in V_i$, where $i\in \{2,3,4\}$.  Clearly, $0\le|\mathcal{B'}|\le 3$.

{\bf CASE 1.}\ \ {\it $|\mathcal{B'}|=0$.} 
Since $V_i\notin \mathcal{B'}$, there is an accessible path $P$ from $V_1$ to $V_i$. 
We move $u$ to $V_5$, move $v$ to $V_1$ and move vertices along the  path $P$. Since $d_G(u)<d_G(v)$, this contradicts~\eqref{min5}.

{\bf CASE 2.}\ \ {\it $|\mathcal{B'}|=1$.} We may let $\mathcal{B'}=\{V_2\}$.  Then $e[V_1, V_2]\ge s$ and $V_3,V_4$ are reachable from $V_1$. This implies $V_3V_2,V_4V_2\not\in E(D)$. Thus $e[V_3\cup V_4,V_2]\ge 2p$. 
In summary, we have the following:
$$e[V_1, V_2]\ge s,\qquad e[V_1, V_5]\ge s,\qquad
e[V_3\cup V_4, V_2]\ge 2p,\qquad 
e[V_2\cup V_3\cup V_4, V_5]\ge 3p 
.$$
\noindent Hence $e[V_1\cup V_3\cup V_4, V_2\cup V_5]+e(G[V_2\cup V_5])\ge(2s+4p)+p=2s+5p$. If $s<p$, then $e(G)\ge2s+5p>7s+1 $, a contradiction. Assume that $p\le s\le 3p-1$. If $G[V_1\cup V_3\cup V_4]$ contains three disjoint independent $p$-sets, then $G[V_1\cup V_2\cup V_3\cup V_4]$ contains four disjoint independent $p$-sets, a contradiction. If $G[V_1\cup V_3\cup V_4]$ contains no three disjoint independent $p$-sets, then by Theorem 1.7, $e(G[V_1\cup V_3\cup V_4])\ge 5(s-p+1)$. Thus $$e(G)=e(G[V_1\cup V_3\cup V_4])+e[V_1\cup V_3\cup V_4, V_2\cup V_5]+e(G[V_2\cup V_5])\ge 5(s-p+1)+2s+5p=7s+5>7s+1,$$ a contradiction.

{\bf CASE 3.}\ \ {\it $|\mathcal{B'}|=2$.} We may let $\mathcal{B'}=\{V_2,V_3\}$. { Then $V_1V_4\in E(D)$. Let $v_1\in V_1$ be movable to $V_4$. If $e[V_4,V_5]=p$, then some vertex $v_5\in V_5$ is a solo neighbor of two vertices $v_4,v_4'\in V_4$. 
{ In this case, we move $v_1$ to $V_4$, $v_4$ to $V_5$ and {$v_5$} to $V_1$.}
 Now $v_4'$ is movable to $V_5$, and so $|\mathcal{B}|\le 3$. This contradicts~\eqref{|B|}. 
 
 Thus $e[V_4,V_5]\ge p+1$. Then by Claim 3.1, $e[V_1\cup V_4, V_2\cup V_3\cup V_5]\ge 3(s+p)+1$} and  $e[V_2\cup V_3, V_5]\ge 2p$.
If there are at most $4$ vertices in $V_1$ that are movable to $V_4$, then $e[V_1, V_4]\ge s-4$. In  this case, { $e(G)\ge 4s+5p-3$}. If $s<p$, then { $4s+5p-3>7s+1$ }, a contradiction. Assume that $p\le s\le 3p-1$. If $G[V_1\cup V_2\cup V_3]$ contains three disjoint independent $p$-sets, then $G[V_1\cup V_2\cup V_3\cup V_4]$ contains four disjoint independent $p$-sets, a contradiction. If $G[V_1\cup V_2\cup V_3]$ contains no three disjoint independent $p$-sets, then by Theorem 1.7, $e(G[V_1\cup V_2\cup V_3])\ge 5(s-p+1)$. 
 By Claim 3.1,  $$e[V_1\cup V_2\cup V_3, V_4\cup V_5]+e(G[V_4\cup V_5])=e[V_1\cup V_2\cup V_3\cup V_4, V_5]+e[V_4,  V_2\cup V_3]+e[V_1, V_4]\ge { 2s+5p-3. }$$ Thus { $$e(G)=e(G[V_1\cup V_2\cup V_3])+e[V_1\cup V_2\cup V_3, V_4\cup V_5]+e(G[V_4\cup V_5])\ge 5(s-p+1)+2s+5p-3=7s+2,$$} a contradiction. 

Assume now that there are at least $5$ vertices that are movable to $V_4$. Since $|\mathcal{B'}|=2$ and $\Delta(G)\le 6$, $\Delta(G[V_1])\le 3$. Then there are  non-adjacent vertices $v_1, v_1'\in V_1$ that are movable to $V_4$.
If $v_2\in V_2$ is a solo neighbor of  vertices $v_4,v_4'\in V_4$, then we move $v_1$ to $V_4$, $v_4$ to $V_2$ and $v_2$ to $V_1$.
Now $v_4'$ is movable to $V_2$. Since $v_1'$ is movable to $V_4$, there is an accessible path $V_1V_4V_2$. Then $|\mathcal{B'}|\le 1$, contradicting~\eqref{min5}.

So, we may assume that no $x\in V_2$  is  a solo neighbor of two vertices in $V_4$. 
 Consider $e[V_4, V_2]+e[V_2, V_3]$. For every vertex $a_4\in V_4$, $a_4$ has at least one neighbor in $V_2$. If $a_4$ has a solo neighbor $a_2\in V_2$, then $a_2$ has neighbors in $V_3$. 
Otherwise let $v_1,v_1'\in V_1$ be two non-adjacent vertices that are movable to $V_4$. We can move $v_1$ to $V_4$, $a_4$ to $V_2$ and $a_2$ to $V_1$.
Now $v_1'$ is movable to $V_4$ and $a_2$ is movable to $V_3$. 
 This implies $V_3\not\in\mathcal{B'}$. So $|B'| \le 1$. This contradicts~\eqref{min5}.

If there are $x$ vertices in $V_4$ having the solo neighbor in $V_2$, then $e[V_4, V_2]+e[V_2, V_3]\ge 2x+2(p-x)=2p$. In summary, we have the following:
$$
e[V_1, V_2\cup V_3\cup V_5]\ge 3s,\qquad
e[V_4, V_2]+e[V_2, V_3]\ge 2p, \qquad
e[V_4, V_3]\ge p,\quad
{ e[V_2\cup V_3\cup V_4, V_5]\ge 3p+1.}
$$

If $s<p$, then { $e(G)\ge 3s+6p+1> 7s+1$}, a contradiction. Assume that $p\le s\le 3p-1$. If $G[V_1\cup V_4]$ contains two disjoint independent $p$-sets, then $G[V_1\cup V_2\cup V_3\cup V_4]$ contains four disjoint independent $p$-sets, a contradiction. If $G[V_1\cup V_4]$ contains no two disjoint independent $p$-sets, by Theorem 1.4, then $e(G[V_1\cup V_4])\ge 3(s-p+1)$. Note that $$\begin{array}{ll}&e[V_1\cup V_4, V_2\cup V_3\cup V_5]+e(G[V_2\cup V_3\cup V_5])\\
=&(e[V_1, V_2]+e[V_1, V_3]+e[V_1, V_5]+e[V_4, V_2]+e[V_4, V_3]+e[V_4, V_5])+(e[V_2, V_3]+e[V_2, V_5]+e[V_3, V_5])\\
\ge& { 3s+6p+1 }.\end{array}$$
Thus { $$e(G)=e(G[V_1\cup V_4])+e[V_1\cup V_4, V_2\cup V_3\cup V_5]+e(G[V_2\cup V_3\cup V_5])\ge 3(s-p+1)+3s+6p+1=6s+3p+4>7s+1,$$} a contradiction. 

{\bf CASE 4.}\ \ {\it $|\mathcal{B'}|=3$.} Clearly, $\mathcal{B'}=\{V_2, V_3, V_4\}$. Then $V_1V_i\notin E(D)$ for $2\le i\le 4$. Since $V_1\in\mathcal{B}$ and $\Delta(G)\le 6$, $\Delta(G[V_1])\le 2$.  So, by Claim 3.4, $1\le s\le 3p-3$ and $e(G[V_1])\ge {s\choose 2}-t_{s,p-1}$.

If for some $2\le i\le 4$, $v_i\in V_i$ is a solo neighbor of two  non-adjacent vertices $v_1,v_1'\in V_1$, then 
we move $v_1$ to $V_i$ and $v_i$ to $V_1$. Now
$v_1'$ is movable to $V_i$. This implies $|\mathcal{B'}|\le 2$,     contradicting~\eqref{min5}.

Thus we may assume that no $w\in V_2\cup V_3\cup V_4$ is  a solo neighbor of two non-adjacent vertices in $V_1$. Consider $e[V_1, V_i]+e[V_i, V_j]$, where $2\le i,j\le 4$ and $i\ne j$. Note that every  $a_1\in V_1$ has at least one neighbor in $V_i$.

If a solo neighbor $a_i\in V_i$ of some $a_1\in V_1$ has no neighbors in $V_j\in \mathcal{B}'-V_i$,
then after moving $a_i$ to $V_1$ and $a_1$ to $V_i$, $V_j\not\in\mathcal{B'}$, since now $a_i$ is movable to it. This contradiction shows that for each  $2\leq i\leq 4$, every $a_i\in V_i$ that is
a solo neighbor of some $a_1\in V_1$ has neighbors in each class in
$\mathcal{B}'-V_i$.

Assume that $V_1$  contains $x$ vertices  that have a solo neighbor in $V_i$ and $s-x$ vertices  that have at least two neighbors in $V_i$. Recall that $\Delta(G[V_1])\le 2$. By Theorem 2.2, there are at least $\lceil\frac{x}{3}\rceil$ vertices in $V_1$ that all have  distinct solo neighbors in $V_i$.  Then for any $j\in \{2,3,4\}-\{i\}$,
\begin{equation}\label{ij}
e[V_1,V_i]+ e[V_i,V_j]\ge x+2(s-x)+\lceil\frac{x}{3}\rceil=2s-\lfloor\frac{2x}{3}\rfloor\ge 2s-\frac{2x}{3}\geq \frac{4s}{3}. 
\end{equation}
Applying~\eqref{ij} for pairs $(i,j)\in \{(2,3),(3,4),(4,2)\}$, we get
$e[V_1,V_2\cup V_3\cup V_4]+e(G[V_2\cup V_3\cup V_4])\geq 3\cdot \frac{4s}{3}=4s$. Together with $e(G[V_1])\ge {s\choose 2}-t_{s,p-1}$ and
$e[V_1\cup V_2\cup V_3\cup V_4,V_5]\geq s+3p$, this yields
$$e(G)\ge 5s+3p+{s\choose 2}-t_{s,p-1}=\left\{\begin{array}{lll}{ 5s+3p>7s+1,} & \mbox{if $0\le s<p$,}\\ 5s+3p+s-(p-1)={6s+2p+1>7s+1, } & \mbox{if $p\le s<2(p-1),$}\\ 5s+3p+2s-3(p-1)={ 7s+3>7s+1, } & \mbox{if $2(p-1)\le s\le 3(p-1)$,}\end{array}\right.$$ a contradiction. \ $\Box$

Consider the following discharging procedure. At the start, each $v\in V_1\cup V_2\cup V_3\cup V_4$ has charge $ch(v)=1$, each $v\in V_5$ has charge $ch(v)=0$ and each $e\in E(G)$ has charge $ch(e)=0$. So, $\sum_{x\in V(G)\cup E(G)}ch(x)=s+3p$. 
Denote $V_0:=V_1\cup V_2\cup V_3\cup V_4$ and
$E_0:=E(G[V_0])$. 
Now we will move the charges between edges and vertices without changing the total sum as follows.

Step 1.\ \ Every vertex $v$ in { $V_0$} gives charge $\frac{1}{10}$ to each edge incident to $v$ that is in { $G[V_0]$}.

Step 2.\ \ Every vertex $v$ in { $V_0$} distributes its current charge equally to its neighbors in $V_5$.

Let the charge of each $x\in V(G)\cup E(G)$ after Step $j$ be denoted by $ch_j(x)$. Then $ch_2(e)= 0$ for every $e\in E(V_0,V_5)$ and $ch_2(v)= 0$ for every $v\in V_0$.  Therefore $$s+3p=\sum_{e\in E_0} ch_2(e)+\sum_{v\in V_5}ch_2(v).$$

{\bf Claim 3.15.}\ \ {\it For each $v\in V_5$, $ch_2(v)\le 3$.}

{\bf Proof of Claim 3.15.}\ \ Recall that $\Delta(G)\le 6$. If $d_G(v)\le 3$, then $ch_2(v)\le 3$. Assume that $4\le d_G(v)\le 6$. Let $N(v)=\{u_1,u_2,\ldots, u_{d_G(v)}\}$. If $v$ is the solo neighbor of $u_i\in N(v)$, then by Claim 3.14, $d_G(u_i)\ge d_G(v)$. Hence $ch_2(v)\le d_G(v)\cdot\left[1-(d_G(v)-1)\cdot\frac{1}{10}\right]=d_G(v)\cdot\frac{11-d_G(v)}{10}$. If $v$ is not the solo neighbor of $u_i\in N(v)$, then $ch_2(v)\le d_G(v)\cdot\frac{1}{2}$. Therefore, $ch_2(v)\le d_G(v)\cdot\max\{\frac{11-d_G(v)}{10}, \frac{1}{2}\}\le 3$. \ $\Box$

By Claim 3.15, $\sum_{v\in V_5}ch_2(v)\le 3(p-1)$. For each edge $e\in E_0$, $ch_2(e)\le\frac{1}{10}+\frac{1}{10}=\frac{1}{5}$. Since $$s+3p=\sum_{e\in E_0} ch_2(e)+\sum_{v\in V_5}ch_2(v)\le\frac{1}{5}e(G[V_1\cup V_2\cup V_3\cup V_4])+3(p-1),$$ we have $|E_0|\ge 5(s+3)$. Thus 
{ \begin{equation}\label{B4}
e(G)=|E_0|+e[V_0,V_5]\geq 5(s+3)+(s+3p)=6s+3p+15>7s+1,
\end{equation} }
a contradiction.  \ $\Box$

\section{ Proof of Theorem 1.8}

\hskip\parindent We prove Theorem 1.8 using Theorem 2.7, first proving Parts 1--4 of the theorem and then Part 5 by induction.
{By $H_n$ we denote an extremal
$n$-vertex graph for the original Tur\' an problem, but mostly consider its complement, $\ol{H}_n$.}

{\bf Proof of Theorem 1.8.}\ \ For $n\ge 4p+1$, we let $s=n-4p+1$ and $\ol{H}_{n}\in\ol{\rm EX}(n, 4K_{p})$. Recall that $K_{s+4}\cup\ol{K}_{4p-5}$ has no four disjoint independent $p$-sets for $2\leq s\leq 5$. Then $\ol{\rm ex}(4p+1,4K_p)\le 15$, $\ol{\rm ex}(4p+2,4K_p)\le 21$, $\ol{\rm ex}(4p+3,4K_p)\le 28$ and $\ol{\rm ex}(4p+4,4K_p)\le 36$. Moreover, since $K_8\cup S_7$ contains no four disjoint independent $3$-sets, $\ol{\rm ex}(4p+4,4K_p)\le e(K_8\cup S_7)=35$ for $p=3$. By Claim 2.1, $\ol{\rm ex}(n,4K_p)\le 7s$ for $4p+5\leq n\leq 8p-8$.

{\em Proof of Part 1.} 
If $\Delta(\ol{H}_{4p+1})\le 6$, then by Theorem 2.7
for $s=2$, $\ol{H}_{4p+1}=K_6\cup\ol{K}_{4p-5}$. If $\Delta(\ol{H}_{4p+1})\ge 7$, let $v\in V(\ol{H}_{4p+1})$ be a vertex with $d_{\ol{H}_{4p+1}}(v)=\Delta(\ol{H}_{4p+1})$ and $G=\ol{H}_{4p+1}\setminus\{v\}$. Then $e(G)\le 15-7=8$ and $|V(G)|=4p$. By Theorem 1.6, 
$\ol{{\rm ex}}(4p, 4K_p)=10$. But $e(G)\le 8$. Thus $G$ contains four disjoint independent $p$-sets, and hence $\ol{H}_{4p+1}$ has four disjoint independent $p$-sets, a contradiction. Thus $e(\ol{H}_{4p+1})=15$ and $\ol{H}_{4p+1}=K_6\cup\ol{K}_{4p-5}$.


{\em Proof of Part 2.} If $\Delta(\ol{H}_{4p+2})\le 6$, then by Theorem 2.7 for $s=3$, $\ol{H}_{4p+2}=K_7\cup\ol{K}_{4p-5}$. If $\Delta(\ol{H}_{4p+2})\ge 7$, let $v\in V(\ol{H}_{4p+2})$ be a vertex with $d_{\ol{H}_{4p+2}}(v)=\Delta(\ol{H}_{4p+2})$ and $G=\ol{H}_{4p+2}\setminus\{v\}$. Then $e(G)\le 21-7=14$ and $|V(G)|=4p+1$. By Part 1, $\ol{{\rm ex}}(4p+1, 4K_p)=15> e(G)$. Thus $G$ contains four disjoint independent $p$-sets, and hence $\ol{H}_{4p+2}$ has four disjoint independent $p$-sets, a contradiction. So $e(\ol{H}_{4p+2})=21$ and $\ol{H}_{4p+2}=K_7\cup\ol{K}_{4p-5}$.


{\em Proof of Part 3.} By Theorem 2.7 for $s=4$, $\Delta(\ol{H}_{4p+3})\geq 7$.
If $\Delta(\ol{H}_{4p+3})=7$, let $v\in V(\ol{H}_{4p+3})$ be a vertex with $d_{\ol{H}_{4p+3}}(v)=\Delta(\ol{H}_{4p+3})=7$,
contained in the most copies of  $K_7$. Let $W=N_{\ol{H}_{4p+3}}(v)=\{v_1,v_2,\ldots, v_7\}$ and consider $G=\ol{H}_{4p+3}\setminus\{v\}$. Note that $e(G)\le 28-7=21$ and $|V(G)|=4p+2$. By Part 2, $G 
=K_7\cup\ol{K}_{4p-5}$. { Let $U=V(K_7)=\{u_1,u_2,\ldots,u_7\}$.}
If $|W\cap U|=0$, then $\ol{H}_{4p+3}=K_7\cup S_7\cup \ol{K}_{4p-12}$. If $4p-12>0$, let $
\{a_1,a_2,\ldots,a_{4p-12}\}=V(G)\setminus (U\cup W)$. 
 Then $\{v,u_1,a_1,a_2,\ldots,a_{p-2}\}$, $\{v_1,v_2,u_2,a_{p-1},a_p,\ldots,a_{2p-5}\}$, $\{v_3,v_4,u_3,a_{2p-4},a_{2p-3},\ldots,a_{3p-8}\}$ and $\{v_5,v_6,v_7,u_4,a_{3p-7},a_{3p-6},\ldots,a_{4p-12}\}$ are $4$ independent sets of size  $p$ in $\ol{H}_{4p+3}$, 
a contradiction. Hence $4p-12=0$, i.e. $p=3$ and $\ol{H}_{4p+3}=K_7\cup S_7$. 
   
If $1\le |N_{\ol{H}_{4p+3}}(v)\cap U|\le 6$, then each vertex in $U\cap N(v)$ has degree $7$ and is contained in more copies of $K_7$ than $v$, contradicting the choice of $v$.
If $|N_{\ol{H}_{4p+3}}(v)\cap V(K_7)|=7$, then $\ol{H}_{4p+3}=K_8\cup\ol{K}_{4p-5}$.

If $\Delta(\ol{H}_{4p+3})\ge 8$, let $v'\in V(\ol{H}_{4p+3})$ be a vertex with $d_{\ol{H}_{4p+3}}(v')=\Delta(\ol{H}_{4p+3})$ and $G'=\ol{H}_{4p+3}\setminus\{v'\}$. Then $e(G')\le 28-8=20$ and $|V(G')|=4p+2$. By Part 2, 
$G'$ contains four disjoint independent $p$-sets, and hence $\ol{H}_{4p+3}$ has four disjoint independent $p$-sets, a contradiction.


{\em Proof of Part 4.} Recall that $\ol{\rm ex}(4p+4,4K_p)\le\left\{\begin{array}{lll} 
35, & \mbox{if $p=3$,}\\
36,  & \mbox{if $p\ge 4$.}\end{array}\right.$ 

If $\Delta(\ol{H}_{4p+4})\le 6$, then by Theorem 2.7 for $s=5$, 
 $\ol{H}_{4p+4}=K_6\cup K_7\cup\ol{K}_{4p-9}$. 

If $\Delta(\ol{H}_{4p+4})\ge 9$, let $v\in V(\ol{H}_{4p+4})$ be a vertex with $d_{\ol{H}_{4p+4}}(v)=\Delta(\ol{H}_{4p+4})$ and $G=\ol{H}_{4p+4}\setminus\{v\}$. Then $e(G)\le\left\{\begin{array}{lll} 
26, & \mbox{if $p=3$,}\\
27,  & \mbox{if $p\ge 4$,}\end{array}\right.$ and $|V(G)|=4p+3$. By Part 3, $\ol{{\rm ex}}(4p+3, 4K_p)=28>e(G)$. Thus $G$ contains four disjoint independent $p$-sets, and hence $\ol{H}_{4p+4}$ has four disjoint independent $p$-sets, a contradiction. 

If $\Delta(\ol{H}_{4p+4})=8$, let $v\in V(\ol{H}_{4p+4})$ be a vertex with $d_{\ol{H}_{4p+4}}(v)=8$.  Consider $G=\ol{H}_{4p+4}\setminus\{v\}$. We have $e(G)\le\left\{\begin{array}{lll} 
27, & \mbox{if $p=3$,}\\
28,  & \mbox{if $p\ge 4$,}\end{array}\right.$ and $|V(G)|=4p+3$. By Part 3, $p\ge 4$ and $G=\ol{H}_{4p+3}=K_8\cup\ol{K}_{4p-5}$. Suppose 
$V(G)=U\cup W$, where $G[U]=K_8$,
$U=\{u_1,\ldots,u_8\}$ and $W=\{w_1,\ldots,w_{4p-5}\}$. Let $y=|N_{U}(v)|$.
If $y<8$, we may assume  $N_{U}(v)=\{u_1,u_2,\ldots,u_{y}\}$ and  $N_{W}(v)=\{w_{4p-4-y'},w_{4p-3-y'}\ldots,w_{4p-5}\}$, where $y'=8-y$. Since $p\ge 4$ and $y'\le 8$, $4p-4-y'>p-2$. Let $b_1,b_2,b_3\in U\setminus\{u_{y+1}\}$, then $\{v, u_{y+1}, w_1, w_2,\ldots,w_{p-2}\}$, $\{b_1, w_{p-1}, w_{p},\ldots,w_{2p-3}\}$, $\{b_2, w_{2p-2}, w_{2p-1},\ldots,w_{3p-4}\}$ and $\{b_3, w_{3p-3}, w_{3p-2},\ldots,w_{4p-5}\}$ are $4$ independent sets with size $p$ in $\ol{H}_{4p+4}$, a contradiction. If $y=8$, then $\ol{H}_{4p+4}=K_9\cup\ol{K}_{4p-5}$.

The last and most difficult case is $\Delta(\ol{H}_{4p+4})=7$.
First, we 
\begin{equation}\label{K7}
\mbox{assume that}\;  \ol{H}_{4p+4}\;   \mbox{contains } \; K_7, \; \mbox{say with vertex set $U=\{u_1,\ldots,u_7\}$.}
\end{equation}
Let $\ol U:=V(\ol{H}_{4p+4})\setminus U$.
If $|E(U,\ol U)|\leq 6$, consider $G=\ol{H}_{4p+4}\setminus U$. Then $|E(G)|\leq 15-|E(U,\ol U)|$. By Theorem 1.2 and Part 1, either $G$ contains four disjoint independent $(p-1)$-sets or $p\geq 4$, $G=K_6\cup \ol{K}_{4p-9}$ and $|E(U,\ol U)|=0$. In the former case, we can add to each of the four disjoint independent $(p-1)$-sets a vertex from $U$ not adjacent to any vertex in it, contradicting the definition of $\ol{H}_{4p+4}$. In the latter, 
$\ol{H}_{4p+4}=K_7\cup K_6\cup \ol{K}_{4p-9}$ satisfies our theorem.

Suppose now $|E(U,\ol U)|=7$. If there is $u_0\in \ol U$ adjacent to all vertices in $U$, then
$\ol{H}_{4p+4}[U+u_0]=K_8$ is a component in $\ol{H}_{4p+4}$, and $\ol{H}_{4p+4}\setminus { \{U+u_0\} }$
does not contain four disjoint independent $(p-1)$-sets. Since $|E(\ol{H}_{4p+4}\setminus { \{U+u_0\} })|<10 $, { by 
Theorem 1.6 for $k=4$,} we have $p-1=2$ and by 
Theorem 1.2, $\ol{H}_{4p+4}\setminus { \{U+u_0\} }=S_7$. Thus $\ol{H}_{4p+4}=K_8\cup S_7$ satisfies our theorem. So, suppose $U$ has at least two neighbors in $\ol U$, say $u_1u'_1,u_2u'_2\in E(\ol{H}_{4p+4})$, where ${ u_1',u'_2}\in \ol U$.
Let $G'$ be obtained from $\ol{H}_{4p+4}\setminus U$ by adding edge $u'_1u'_2$. 
Then $G'$ has $4p-3=4(p-1)+1$ vertices and at most $36-28+1=9$ edges. By Theorem 1.2 and Part 1, $G'$ has four disjoint independent $(p-1)$-sets, say $W_1,\ldots,W_4$. As in the proof of~\eqref{10}, 
if there are four vertices
$u_1,\ldots,u_4\in U$ such that for $j=1,2,3,4$,  $u_j$ has no neighbors in $W_j$, then the sets
$W'_j=W_j+u_j$ form four disjoint independent $p$-sets, a contradiction. Otherwise,  
since  each $u\in U$ has at most one neighbor outside of $U$, by Hall's Theorem, 
there is $j\in [4]$ such that $N(W_j)\supseteq U$.
But then $u_1',u'_2\in W_j$, contradicting the fact that $W_j$ is independent.
So, all cases when~\eqref{K7}  holds are considered. 

Suppose now
\begin{equation}\label{K72}
  \ol{H}_{4p+4}\;   \mbox{does not contain } \; K_7.
\end{equation}
 Let $X=\{v_1,\ldots,v_x\}$ be a smallest inclusion maximal independent set of degree $7$ vertices, and let $G=\ol{H}_{4p+4}\setminus X$. By the maximality of $X$,
\begin{equation}\label{K71}
   \mbox{$\Delta(G)\leq 6$.}
\end{equation}
If $x=1$, then $|V(G)|=4p+3$ and $e(G)\le\left\{\begin{array}{lll} 
28, & \mbox{if $p=3$,}\\
29,  & \mbox{if $p\ge 4$.}\end{array}\right.$ So,~\eqref{K72} together with~\eqref{K71} contradict  Theorem 2.7 for $s=4$. 

If $x=2$, then $|V(G)|=4p+2$ and $e(G)\le 22
.$ By Theorem 2.7(c) for $s=3$, $K_7\subseteq G$, a contradiction to~\eqref{K72}.

If $x=3$, then $|V(G)|=4p+1$ and $e(G)\le\left\{\begin{array}{lll} 
14, & \mbox{if $p=3$,}\\
15,  & \mbox{if $p\ge 4$.}\end{array}\right.$ 

Then by Part 1 of our theorem, $p\geq 4$ and
 $G=K_6\cup\ol{K}_{4p-5}$. Suppose the vertex set of $K_6$ is $W=\{w_1,\ldots,w_6\}$ 
 and
 $V(G)\setminus (X\cup W)=Y=\{y_1,\ldots,y_{4p-5}\}$. By~\eqref{K72}, we may assume $v_1w_6\notin E(\ol{H}_{4p+4})$. Since $p\geq 4$, $3p-3>7$, so we may assume that $v_1$ is not adjacent to $y_i$ for $3p-2\leq i\leq 4p-5$. We form four independent sets as follows:
 for $1\leq j\leq 3$, we let $I_j=\{y_{(j-1)(p-1)+1}, y_{(j-1)(p-1)+2},\ldots,y_{j(p-1)}\}\cup \{w_j\}$ and
$$I_4=\{y_{3p-2}, y_{3p-1},\ldots,y_{4p-5}\}\cup \{v_1,w_6\}.
 $$
These four sets contradict the choice of $\ol{H}_{4p+4}$.
 
If $x\ge 4$, let $G'=\ol{H}_{4p+4}\setminus\{v_1,v_2,\ldots, v_4\}$. Then $|V(G')|=4p$ and $e(G')\le36-28=8$. By Theorem 1.6, 
$\ol{{\rm ex}}(4p, 4K_p)=10$. Then $G'$ contains four disjoint independent $p$-sets. Thus $\ol{H}_{4p+4}$ has four disjoint independent $p$-sets, a contradiction. This proves Part 4.

{\em Proof of Part 5.} We need a couple of lemmas.

{\bf Lemma 4.1.} \cite{Zhang2024DisjointCliques}\ \ {\it Let $G$ be a $kK_p$-free graph on $n\ge kp$ vertices with $p\ge 3$ and $k\ge 2$. Then $\Delta(\ol{G})\ge\lfloor\frac{n-k}{p-1}\rfloor$ and $\delta(G)\le n-1-\lfloor\frac{n-k}{p-1}\rfloor$.}

{\bf Lemma 4.2.}\ \ {\it Let $p\ge 3$ and $k\ge 2$. If there is an integer $n_1\ge kp$ such that $EX(n_1, kK_p)=\left\{K_{k-1}\vee T_{n_1-k+1,p-1}\right\}$, then $EX(n, kK_p)=\left\{K_{k-1}\vee T_{n-k+1,p-1}\right\}$ for all $n\ge n_1$.}

{\bf Proof.} \ \ We use induction on $n$. Let $G\in EX(n, kK_p)$ for $n>n_1$. By Lemma 4.1, there is a vertex $v\in V(G)$ with $d(v)=\delta(G)\le n-1-\lfloor\frac{n-k}{p-1}\rfloor$. By Lemma 2.5, ${\rm ex}(n, kK_p)=e(K_{k-1}\vee T_{n-k+1,p-1})$ for every $n\ge n_1$. Then $$e(K_{k-1}\vee T_{n-k+1,p-1})=e(G)=e(G-v)+d(v)\le e(K_{k-1}\vee T_{n-k,p-1})+n-1-\lfloor\frac{n-k}{p-1}\rfloor=e(K_{k-1}\vee T_{n-k+1,p-1}).$$ Thus $d(v)=n-1-\lfloor\frac{n-k}{p-1}\rfloor$ and $e(G-v)=e(K_{k-1}\vee T_{n-k,p-1})$. By induction hypothesis, $G-v=K_{k-1}\vee T_{n-k,p-1}$. Then $V(G-v)$ can be divided into $p$ classes $V',V_1,V_2,\ldots,V_{p-1}$ such that $|V'|=k-1$ and $\lfloor\frac{n-k}{p-1}\rfloor\le |V_i|\le \lceil\frac{n-k}{p-1}\rceil$ for $1\le i\le p-1$. Note that $(G-v)[V']=K_{k-1}$ and $(G-v)[V_1\cup V_2\cup \cdots\cup V_{p-1}]=T_{n-k,p-1}$. If $N(v)\cap V_i\not=\emptyset$ for $1\le i\le p-1$, let $v_i\in N(v)\cap V_i$. Then $G[v,v_1,v_2,\ldots,v_{p-1}]=K_p$. Note that $G-\{v,v_1,v_2,\ldots,v_{p-1}\}=K_{k-1}\vee T_{n-p-k+1,p-1}$. Since $n\ge kp$, $G-\{v,v_1,v_2,\ldots,v_{p-1}\}$ contains a copy of $(k-1)K_p$. Then $G$ contains a copy of $kK_p$, a contradiction. Thus there is a $1\le i_0\le p-1$ such that $N(v)\cap V_{i_0}=\emptyset$. Since $$n-1-\lfloor\frac{n-k}{p-1}\rfloor=d(v)\le n-1-|V_{i_0}|\le n-1-\lfloor\frac{n-k}{p-1}\rfloor,$$ $N(v)=V(G-v)\setminus V_{i_0}$, $|V_{i_0}|=\lfloor\frac{n-k}{p-1}\rfloor$ and $G=K_{k-1}\vee T_{n-k+1,p-1}$. \ $\Box$

By Lemma 4.2, we have the following corollary.

{\bf Corollary 4.3.}\ \ {\it If there is an integer $n_1\ge kp$ such that $\ol{EX}(n_1, kK_p)=\left\{\ol{K}_{k-1}\cup\ol{T}_{n_1-k+1,p-1}\right\}$, then $\ol{EX}(n, kK_p)=\left\{\ol{K}_{k-1}\cup\ol{T}_{n-k+1,p-1}\right\}$ for all $n\ge n_1$.
\ $\Box$}

First, we show that 
$\ol{\rm ex}(n,4K_p)=7(n-4p+1)$ for $4p+5\le n\le 7p-2$. In other words, we prove:
\begin{equation}\label{43}
    \mbox{\em For $6\le s\le 3p-1$ and $n=4p+s-1$, $e(\ol{H}_n)=7s$.}
\end{equation}
For this, recall  the definitions of $\JJ_n(p)$ and $\JJ'_n(p)$  given in Section 1. Let $s=n-4p+1$. By Claim 2.1, $e(\ol{H}_n)\le7s$. To prove the lower bound on $e(\ol{H}_n)$, we use induction on $s$.

{ The base case is $s=6$, i.e., $n=4p+5$. Then $e(\ol{H}_{4p+5})\le 42$. If $\Delta(\ol{H}_{4p+5})\le 6$, then by Theorem 2.7 for $s=6$, $e(\ol{H}_{4p+5})=42$ and $\ol{H}_{4p+5}=2K_7\cup\ol{K}_{4p-9}$. If $\Delta(\ol{H}_{4p+5})\ge 7$, let $v\in V(\ol{H}_{4p+5})$ be a vertex with $d_{\ol{H}_{4p+5}}(v)=\Delta(\ol{H}_{4p+5})$ and $G=\ol{H}_{4p+5}\setminus\{v\}$. Then $e(G)\le 42-7=35$ and $|V(G)|=4p+4$. If $e(G)=35$, then $d_{\ol{H}_{4p+5}}(v)=\Delta(\ol{H}_{4p+5})=7$. By Part 4, $p=3$ and $G=K_8\cup S_7$. Suppose that the set of vertices of $G$ is $\{u_1,u_2,\ldots,u_8,w_1,w_2,\ldots,w_8\}$, where $G[\{u_1,u_2,\ldots,u_8\}]=K_8$, $G[\{w_1,w_2,\ldots,w_8\}]=S_7$ and $d_G(w_1)=7$. Thus $N_{\ol{H}_{4p+5}}(v)=\{w_2,\ldots,w_8\}$. Note that $\{v,u_1,w_1\}$, $\{u_2,w_2,w_3\}$, $\{u_3,w_4,w_5\}$ and $\{u_4,w_6,w_7\}$ are $4$ independent sets of size $3$. Thus $\ol{H}_{17}$ contains four disjoint independent $3$-sets, a contradiction. If $e(G)\le 34$, then by Part 4, $G$ contains four disjoint independent $p$-sets, and hence $\ol{H}_{4p+5}$ has four disjoint independent $p$-sets, a contradiction. }

Assume $s\ge 7$. 
If $\Delta(\ol{H}_n)\le 6$, then by Theorem 2.7, $e(\ol{H}_n)\geq 7s$. Suppose $\Delta(\ol{H}_n)\geq 7$ and $e(\ol{H}_n)\leq 7s-1$. Choose $v\in V(\ol{H}_n)$ with $d_{\ol{H}_n}(v)\geq 7$. Then $e(\ol{H}_n\setminus\{v\})\le (7s-1)-7=7s-8$. Since $\ol{H}_n\setminus\{v\}$ contains no four disjoint independent $p$-sets, this contradicts the
induction hypothesis. 

By~\eqref{43}, $\ol{\rm ex}(n,4K_p)=7s=7(n-4p+1)$ for $4p+5\le n\le 7p-2$. Note that $\ol{\rm ex}(7p-2,4K_p)=21p-7$ and $e(\ol{K}_3\cup\ol{T}_{7p-5,p-1})=21(p-3)+56=21p-7$. Thus ${\rm ex}(7p-2,4K_p)=e(K_3\vee T_{7p-5,p-1}).$ By Lemma 2.5, ${\rm ex}(n,4K_p)=e(K_3\vee T_{n-3,p-1})=3+3(n-3)+t_{n-3,p-1}$ for $n\ge 7p-1$. Thus $$\ol{\rm ex}(n,4K_p)=\left\{\begin{array}{lll} 
7(n-4p+1), & \mbox{if $4p+5\le n\le 7p-2$,}\\
\ol{t}_{n-3,p-1},  & \mbox{if $n\ge 7p-1$.}\end{array}\right.$$ Next we determine $\ol{\rm EX}(n,4K_p)$. By Claim 2.1, $\JJ_n(p)\subseteq\ol{\rm EX}(n,4K_p)$ for $4p+5\le n\le 7p-7$ and $\JJ'_n(p)\subseteq\ol{\rm EX}(n,4K_p)$ for $7p-6\le n\le 8p-8$. To prove that $\ol{\rm EX}(n,4K_p)\subseteq\JJ_n(p)$ and $\ol{\rm EX}(n,4K_p)\subseteq\JJ'_n(p)$, we use induction on $s$.

{\bf Claim 4.4.} \ \ {\it For $4p+6\le n\le 8p-7$, $\ol{H}_n\in\JJ_n(p)$ or $\Delta(\ol{H}_n)=7$.}

{\bf Proof of Claim 4.4.}\ \ Assume $4p+6\le n\le 7p-2$. If $\Delta(\ol{H}_n)\le 6$, then by~\eqref{43} and Theorem 2.7, $s$ is divisible by $3$ and $\ol{H}_n=\frac{s}{3}K_7\cup\ol{K}_{n-\frac{7s}{3}}\in\JJ_n(p)$. If $\Delta(\ol{H}_n)\ge 8$, let $v\in V(\ol{H}_n)$ be a vertex with $d_{\ol{H}_n}(v)=\Delta(\ol{H}_n)$ and $G=\ol{H}_n\setminus\{v\}$. Then $e(G)\le 7s-8$ and $|V(G)|=4p+s-2$. Since $e(\ol{H}_{n-1})=7s-7$ and $e(G)<e(\ol{H}_{n-1})$, $G$ contains four disjoint independent $p$-sets, and hence $\ol{H}_n$ has four disjoint independent $p$-sets, a contradiction. Thus $\Delta(\ol{H}_n)=7$.

If $7p-2\ge 8p-7$, then we are done. If $7p-2<8p-7$, then by Lemma 4.1, $\Delta(\ol{H}_n)\ge\lfloor\frac{n-4}{p-1}\rfloor\ge 7$ for $7p-1\le n\le 8p-7$. If $\Delta(\ol{H}_n)\ge 8$, then we can delete a vertex with maximum degree so that the number of edges of the resulting graph is less than $e(\ol{H}_{n-1})$. Thus the resulting graph has four disjoint independent $p$-sets, a contradiction. Thus, $\Delta(\ol{H}_n)=7$. \ $\Box$

By Claim 4.4, we may suppose $\Delta(H_n)=7$. Let $v\in V(\ol{H}_n)$ be a vertex with $d_{\ol{H}_n}(v)=\Delta(\ol{H}_n)=7$, contained in the most copies of $K_7$. Consider $G=\ol{H}_n\setminus\{v\}$.

{\bf Claim 4.5.} \ \ {\it (a)\ If $4p+6\le n\le 7p-7$ and $G\in\JJ_{n-1}(p)$, then $\ol{H}_n\in\JJ_{n}(p)$;\\
(b)\ If $n=7p-6$ and $G\in\JJ_{7p-7}(p)$, then $\ol{H}_n\in\JJ'_{n}(p)$;\\
(c)\ If $7p-5\le n\le 8p-8$ and $G\in\JJ'_{n-1}(p)$, then $\ol{H}_n\in\JJ'_{n}(p)$.
}

{\bf Proof of Claim 4.5.}\ \ If $G\in\JJ_{n-1}(p)$, let $G=x_1K_7\cup y_1K_8\cup\ol{K}_{n-1-7x_1-8y_1}$. {Since $\Delta(H_n)=7$, $v$ has no neighbors in the copies of $K_8$.
If $v$ is adjacent to some vertex $u$ in a copy $F$ of $K_7$ in $G$, but is not adjacent to some other vertex $u'\in V(F)$, then $u$ is in more copies of $K_7$ than
$v$. This contradicts the choice of $v$. Thus}
all neighbors of $v$ are either contained in the same copy of $K_7$ of $G$ or are isolated in $G$.
 If the former holds, then $\ol{H}_n=(x_1-1)K_7\cup (y_1+1)K_8\cup\ol{K}_{n-1-7x_1-8y_1}\in\JJ_{n}(p)$. In the latter case, $n-1-7x_1-8y_1\ge 7$. If $n-1-7x_1-8y_1>7$, then $\ol{H}_n=x_1K_7\cup y_1K_8\cup S_7\cup\ol{K}_{n-7x_1-8y_1-8}$ and $n-7x_1-8y_1-8\ge 1$. Thus $S_7\cup\ol{K}_{n-7x_1-8y_1-8}$ has an equitable $4$-coloring. Hence $x_1K_7\cup y_1K_8\cup S_7\cup\ol{K}_{n-7x_1-8y_1-8}$ contains four disjoint independent $\left(x_1+y_1+\frac{n-7x_1-8y_1}{4}\right)$-sets. Since $x_1+y_1+\frac{n-7x_1-8y_1}{4}=\frac{n-3x_1-4y_1}{4}=p$, $\ol{H}_n$ has four disjoint independent $p$-sets, a contradiction. If $n-1-7x_1-8y_1=7$, then $x_1=8p-8-n$ and $y_1=n-7p+6$. Thus $7p-6\le n\le 8p-8$ and $\ol{H}_n=(8p-8-n)K_7\cup (n-7p+6)K_8\cup S_7\in\JJ'_{n}(p)$. 

If $G=(8p-7-n)K_7\cup (n-7p+5)K_8\cup S_7$, then again by the choice of $v$, all neighbors of $v$ are either contained in the same $K_7$ of $G$ or in $S_7$ of $G$. If the former holds, then $\ol{H}_n=(8p-8-n)K_7\cup (n-7p+6)K_8\cup S_7\in\JJ'_{n}(p)$. In the latter case, $\ol{H}_n=(8p-7-n)K_7\cup (n-7p+5)K_8\cup K_{2,7}$ contains four disjoint independent $p$-sets, a contradiction.  \ $\Box$

{\bf Claim 4.6.} \  {\it 
(a) $\ol{\rm EX}(n,4K_p)=\JJ_n(p)$ if $4p+5\le n\le 7p-7$; $\,$ (b)
 $\ol{\rm EX}(n,4K_p)=\JJ'_n(p)$ if $7p-6\le n\le 8p-8$.
}

{\bf Proof of Claim 4.6.}\ \ If $n=4p+5$, then $s=6$. Recall that $\ol{H}_{4p+5}=2K_7\cup\ol{K}_{4p-9}\in \JJ_{4p+5}(p)$. If $4p+6\le n\le 7p-7$, then $7\le s\le 3p-6$. Note that $e(G)=7s-7$ and $|V(G)|=n-1$. By the induction hypothesis, $G\in\ol{\rm EX}(n-1,4K_p)\subseteq\JJ_{n-1}(p)$. By Claim 4.5(a), $\ol{\rm EX}(n,4K_p)=\JJ_n(p)$ for $4p+5\le n\le 7p-7$.
This proves (a).

If $n=7p-6$, then $s=3p-5$ and $|V(G)|=7p-7$. By~\eqref{43}, $e(G)=7s-7$. By (a), $G\in \JJ_{7p-7}(p)$. By Claim 4.5(b), $\ol{H}_{7p-6}\in\JJ'_{7p-6}(p)$. If $7p-5\le n\le 8p-8$, then $3p-4\le s\le 4p-7$. Note that $e(G)=7s-7$ and $|V(G)|=n-1$. By the induction hypothesis, $G\in\ol{\rm EX}(n-1,4K_p)\subseteq\JJ'_{n-1}(p)$. By Claim 4.5(c), $\ol{\rm EX}(n,4K_p)=\JJ'_n(p)$ for $7p-6\le n\le 8p-8$.\ $\Box$



If $n=8p-7$, then by Claim 4.4, $\ol{H}_{8p-7}\in\JJ_{8p-7}(p)$ or $\Delta(\ol{H}_{8p-7})=7$. If the former holds, then by the definition of $\JJ_n(p)$, $\JJ_{8p-7}(p)=\{2K_7\cup (p-3)K_8\cup \ol{K}_3\}$. Thus $\ol{H}_{8p-7}=2K_7\cup (p-3)K_8\cup \ol{K}_3$. In the latter case, $\Delta(H_{8p-7})=7$. Let $v'\in V(\ol{H}_{8p-7})$ be a vertex with $d_{\ol{H}_{8p-7}}(v')=\Delta(\ol{H}_{8p-7})=7$, contained in the most copies of $K_7$. Consider $G'=\ol{H}_{8p-7}\setminus\{v'\}$. Note that $e(G')=7s-7$ and $|V(G')|=8p-8$. By Claim 4.6, $G'\in \JJ'_{8p-8}(p)$. By definition of $\JJ'_n(p)$, $\JJ'_{8p-8}(p)=\{3K_7\cup (p-4)K_8\cup \ol{K}_3, (p-2)K_8\cup S_7\}$. By the choice of $v'$, $\ol{H}_{8p-7}\in\{2K_7\cup (p-3)K_8\cup\ol{K}_3, (p-2)K_8\cup K_{2,7}\}$. Since $(p-2)K_8\cup K_{2,7}$ contains four disjoint independent $p$-sets, $\ol{H}_{8p-7}=2K_7\cup (p-3)K_8\cup\ol{K}_3$. Note that $\ol{K}_{3}\cup\ol{T}_{8p-10,p-1}=2K_7\cup (p-3)K_8\cup\ol{K}_3$. Thus $\ol{\rm EX}(8p-7,4K_p)=\{\ol{K}_{3}\cup\ol{T}_{8p-10,p-1}\}$. By Corollary 4.3 for $k=4$, $\ol{\rm EX}(n,4K_p)=\{\ol{K}_{3}\cup\ol{T}_{n-3,p-1}\}$ for $n\ge 8p-7$. This completes the proof of Theorem 1.8. \ $\Box$

\section{Concluding remarks}

\hskip\parindent 1. We think that the following analog of Theorem 2.7 holds.

{\bf Conjecture 5.1.}\ \ {\it Let $p\ge 3$, $k\ge 2$ and $G$ be a graph on $n=kp-1+s$, where $1\le s\le (k-1)p-1$. If $|E(G)|\le (2k-1)s$ and $\Delta(G)\le 2k-2$, then $G$ contains $k$ disjoint independent $p$-sets or $(k-1)|s$ and $G=\frac{s}{k-1}K_{2k-1}\cup\ol{K}_{n-\frac{(2k-1)s}{k-1}}$ with $|E(G)|=(2k-1)s$.}

2. It is likely that our main result can be extended to graphs with no $k$ disjoint $K_p$ as follows.

{\bf Conjecture 5.2.}\ \ {\it Let $p\ge 3$, $k\ge 4$ and $(k-1)p-k^2+3k-3\ge 0$. Then $${\rm ex}(n,kK_p)=\left\{\begin{array}{lll} 
{n\choose2}-(2k-1)(n-kp+1), & \mbox{if $kp+k^2-3k+1\le n\le (2k-1)p-2$,}\\
e(K_{k-1}\vee T_{n-k+1,p-1}),  & \mbox{if $n\ge (2k-1)p-1$.}\end{array}\right.$$}

\paragraph{Acknowledgment.} We thank  Zimu Xiang for fruitful discussions and both  referees for the helpful comments.

\printbibliography

\end{document}